\RequirePackage{fix-cm}
\pdfoutput=1
\documentclass[leqno, fleqn, centertags, 12pt]{article}

\usepackage{latexsym}
\usepackage{amsmath}

\usepackage[T1]{fontenc}
\usepackage[upright]{fourier}

\usepackage{fullpage}

\newcommand{\dis}{\displaystyle}

\newcommand{\rr}{\mathbf{R}}

\newcommand{\toto}{\longrightarrow}

\newcommand{\rl}{ \mathbf{R}^{\dff\mathcal{L}} }
\newcommand{\card}{\mathop{\mbox{\textup{Card}\fff}}}
\newcommand{\dbe}{de~Bruijn--Erd\"{o}s\ }

\def\sss{\hspace{0.05em}\ }
\def\dss{\hspace{0.1em}\ }
\def\trs{\hspace{0.15em}\ }
\def\qss{\hspace{0.2em}\ }
\def\pss{\hspace{0.3em}\ }
\def\oss{\hspace{0.4em}\ }

\def\halfff{\hspace*{0.025em}}
\def\fff{\hspace*{0.05em}}
\def\dff{\hspace*{0.1em}}
\def\trf{\hspace*{0.15em}}
\def\qff{\hspace*{0.2em}}
\def\pff{\hspace*{0.3em}}
\def\off{\hspace*{0.4em}}

\def\ttff{{\hspace*{-0.05em}--\hspace*{0.15em}}}

\newcommand{\nsp}{\hspace*{-0.1em}}
\newcommand{\nnsp}{\hspace*{-0.15em}}

\newcommand{\dnsp}{\hspace*{-0.2em}}

\makeatletter
\renewcommand{\@makefntext}[1]{\vspace*{0.5ex}
\parindent=0em
\hspace*{-0.4em}
\hbox to 0.4em{\hss\@makefnmark}\hspace*{0.4em}{#1}
}
\makeatother

\newcommand{\proof}{\vspace{0.5\bigskipamount}{\textbf{{\emph{Proof}.}}\hspace*{0.7em}}}

\newcommand{\eproof}{ $\blacksquare$}

\newcounter{mysectionnumber}
\newcommand{\mysection}[2]{
\setcounter{equation}{0}\refstepcounter{mysectionnumber}
\section*{ \textnormal{{{\themysectionnumber}.}\oss {#1}}}\label{#2}}

\newcounter{myparnum}[mysectionnumber]
\newcommand{\myuppar}[1]{\vspace{\medskipamount}\textbf{#1}\hspace*{0.5em}}

\numberwithin{equation}{section}

\newcommand{\myit}[1]{\textbf{\textit{#1}}\hspace{0.0em}}

\newcommand{\mynonumsection}[1]{
\vspace{-0.0ex}
\section*{{}\hspace*{0.00em}$\phantom{1.}$\textnormal{{#1}}}}

\title{The\qss de Bruijn--Erd\"{o}s--Hanani\qss theorem}
\author{Nikolai\dss V.\dss Ivanov}
\date{}

\begin{document}

\setlength{\baselineskip}{12pt plus 0pt minus 0pt}
\setlength{\parskip}{12pt plus 0pt minus 0pt}
\setlength{\abovedisplayskip}{12pt plus 0pt minus 0pt}
\setlength{\belowdisplayskip}{12pt plus 0pt minus 0pt}

\newskip\smallskipamount \smallskipamount=3pt plus 0pt minus 0pt
\newskip\medskipamount   \medskipamount  =6pt plus 0pt minus 0pt
\newskip\bigskipamount   \bigskipamount =12pt plus 0pt minus 0pt

\maketitle

\vspace*{10ex}

{\renewcommand{\baselinestretch}{1}
\selectfont

\myit{\hspace*{0em}\large Contents}\vspace*{1.25ex} \\ 
\myit{\phantom{1}1.}\hspace*{0.5em} N.\dss Bourbaki and the\qss
\emph{``Kvant''}\qss magazine \hspace*{0.5em} \vspace*{0.25ex}\\
\myit{\phantom{1}2.}\hspace*{0.5em} A solution of the\dss N.\dss Bourbaki exercise \hspace*{0.5em} \vspace*{0.25ex}\\
\myit{\phantom{1}3.}\hspace*{0.5em} The de Bruijn--Erd\"{o}s proof  \hspace*{0.5em} \vspace*{0.25ex}\\
\myit{\phantom{1}4.}\hspace*{0.5em} From de Bruijn--Erd\"{o}s to systems of distinct representatives  \hspace*{0.5em} \vspace*{0.25ex}\\
\myit{\phantom{1}5.}\hspace*{0.5em} Linear algebra and 
the inequality\qss $m\qff \geqslant\qff n$  \hspace*{0.5em} \vspace*{0.25ex}\\
\myit{\phantom{1}6.}\hspace*{0.5em} Hanani's theorem  \hspace*{0.5em} \vspace*{0.25ex}\\
\myit{\phantom{1}7.}\hspace*{0.5em} Another proof\dss of\dss Hanani's theorem \hspace*{0.5em} \vspace*{0.25ex}\\
\myit{\phantom{1}8.}\hspace*{0.5em} All the de Bruijn--Erd\"{o}s inequalities  \hspace*{0.5em} \vspace*{1.5ex}\\
\myit{References}\hspace*{0.5em}  \hspace*{0.5em}  

}

\footnotetext{\hspace*{-0.65em}\copyright\ 
Nikolai\dss V.\qss Ivanov,\oss 2017.\oss 
Neither the work reported in this paper,\qss 
nor its preparation were 
supported by any governmental 
or non-governmental agency,\qss 
foundation,\qss 
or institution.}

\vspace*{2ex}
\mynonumsection{Preface}

\vspace*{6pt}
{\small 
The present paper is devoted to a somewhat idiosyncratic account of the theorem of 
de Bruijn--Erd\"{o}s and Hanani
from the combinatorics of finite geometries and its various proofs.\oss
Among the proofs discussed are the original proofs by de Bruijn--Erd\"{o}s
and by Hanani\qss (the latter seems to be largely forgotten,\oss
being published in a hard to access journal)
and few others.\oss
Each of these proofs sheds new light on the theorem,\oss
illustrating the maxim that proofs are more important than the theorems proved.\oss
Some proofs and arguments in this paper seem to be new.\oss
I\dss explain how one of the proofs was discovered,\oss
and how another one could have been discovered.\oss
See Sections\qss \ref{reps}\qss and\qss \ref{all}.

I am grateful to\dss F.\dss Petrov for stimulating correspondence
and to\dss M.\dss Prokhorova for careful reading  
of this paper\halfff,\qss
numerous suggestions,\qss 
and providing me with copies of\trs H.\dss Hanani's papers\qss \cite{h1,h2}.

}

\renewcommand{\baselinestretch}{1}
\selectfont

\mysection{N.\qss Bourbaki\qss and\qss the\qss \emph{``Kvant''}\qss magazine}{kvant}

\vspace*{6pt}
\myuppar{Problem.} 
\emph{Let\qss $E$\dss be a set of $n$ elements.\oss
Suppose that $m$ different subsets of\qss $E$
(not equal to\dss $E$ itself\dff)\qss are selected in such a way
that for every two elements of\qss $E$\dss there is exactly one selected subset
containing both these elements.\oss
Prove that\dss $m\qff \geqslant\qff n$\nnsp.\oss}

\emph{When an equality is possible\fff?}

\vspace{6pt}
In\qss 1970\qss this problem was included 
as the Problem\qss M5\qss
in the very first issue of the Soviet\qss \emph{``Kvant''}\qss
magazine and attributed to\dss N.\qss Bourbaki\qss
\cite{b-70}.\qss
The intended audience of the\qss \emph{``Kvant''}\qss magazine\qss
(its name means\qss \emph{``Quantum''}\trf)\qss
was the school students in the USSR of the last two-three grades.\qss
Nowadays the audacity of the editorial board inspires awe:\qss Problem\dss M5\dss was offered
to this audience exactly as it is cited above,\qss as an abstract problem about
finite sets without any motivation and any hints.\qss
The readers were expected to be interested in 
this problem and to appreciate its beauty without any crutches.\qss

In\qss 1970\qss I was among the intended audience of\qss \emph{``Kvant''},\qss
but I was more interested in the foundations of mathematics and in the set theory
than in the combinatorics of finite sets.\qss
I easily found this problem in the Russian translation\qss \cite{b-56}\qss
of the\qss \emph{``Th\'{e}orie des ensembles''}\qss by\qss N.\qss Bourbaki.\qss
It turned out to be the Exercise\qss 12\qss 
to the section\qss \emph{``Calcul sur les entiers''}.\oss
In all editions this exercise is marked as
one of the most difficult\halfff.

The editors of\qss \emph{``Kvant''}\qss were faithful to\dss N.\dss Bourbaki
in not offering any motivation.\qss
But\halfff,\qss in contrast with\qss \emph{``Kvant''},\qss N.\dss Bourbaki split
the result into few steps,\qss offered a hint to the key one,\qss and
stated the expected result in the case of the equality.\qss
The first two steps were rather easy,\qss
but the hint to key step turned out to be incomprehensible for me.\qss

According to the authors of the solution\qss \cite{q-solution}\qss 
published in\qss \emph{``Kvant''}\pss
a few months later\halfff,\qss
they followed\qss 
\emph{``the hints of the author of the problem\fff,\qss N. Bourbaki\dss himself\trf''}\qss
and referred to\qss \cite{b-56}.\qss
The habit of\qss N. Bourbaki to include in his tract recent results
without attribution as exercises is well known,\pss
and was well known in\trs 1970\qss in the Soviet Union.\qss
But it seems that neither the editors of the\qss \emph{``Kvant''}\qss
magazine,\oss nor the authors of the solution\qss \cite{q-solution}\qss
were aware that this result is due to\dss 
N.G.\qss de Bruijn and\dss P.\qss Erd\"{o}s\qss \cite{db-e}\qss
and\dss H.\qss Hanani\qss \cite{h1,h2}.\oss

Neither was I before by an accident I returned to this problem in 2016.\qss 
By this time I was able to immediately recognize that this exercise from\qss \cite{b-56}\qss
is about points and lines in a geometry,\qss
and this realization quickly lead me to the de Bruijn--Erd\"{o}s paper\qss \cite{db-e}.\qss
The exercise turned out to be a quite faithful summary of the de Bruijn--Erd\"{o}s proof\halfff,\qss
and the key part of the proof\halfff,\qss
summarized by\dss N.\qss Bourbaki as a hint\halfff,\qss
turned out to be nearly as obscure as this hint\halfff.
Here is my translation of this exercise
based on the reprint\qss \cite{b-06}\qss
of the\trs 1970\qss edition\qss
(where it appears as the Exercise\qss 14\qss to\qss \S\qss 5).\qss
It is slightly different from the translation in\qss \cite{b-04}.\qss

\myuppar{Exercise.}
\emph{Let\qss $E$\dss be a finite set of $n$ elements,\oss 
$(a_j)_{1\dff \leqslant\dff  j\dff \leqslant\dff  n}$\qss          
be the sequence of elements of $E$ arranged in some order\halfff,\oss 
$(A_i)_{1\dff \leqslant\dff  i\dff \leqslant\dff  m}$\qss     
be a sequence of parts of\qss $E$\nnsp.}

(a)\oss \emph{For each index $j$\nnsp,\oss 
let $k_j$ be the number of indices $i$ such that 
$a_j\qff \in\qff A_i$\nsp;\oss  
for each index $i$ let\oss 
$s_i\qff =\qff \card\dff (\fff A_i\fff)$\dnsp.\qff\oss     
Show that}
\[
\quad
\sum_{j\qff =\qff 1}^n\qff k_j
\off =\off
\sum_{i\qff =\qff 1}^m\qff s_i\off.
\]

\vspace{-9pt}
({\halfff}b)\oss \emph{Suppose that
for each subset\qss $\{\dff x\fff,\pff y\dff\}$\qss
of two elements of\qss $E$\nnsp,\oss
there exists one and only one index $i$ such that $x$ and $y$ are contained in $A_i$\nnsp.\oss
Show that if\qss 
$a_j\qff \not\in\qff A_i$\nnsp,\oss
then\qss $s_i\qff \leqslant\qff k_j$\nnsp.\oss}

(c)\oss \emph{Under the assumptions of\oss \textup{({\halfff}b)},\oss 
show that\qss $m\qff \geqslant\qff n$\nnsp.\oss} 
({\fff}Let $k_n$ be the least of the numbers $k_j$\nnsp;\oss 
show that one may assume that\halfff,\oss whenever\qss 
$i\qff \leqslant\qff k_n\dff,\off j\qff \leqslant\qff k_n$\qss 
and\qss $i\qff \neq\qff j$\nnsp,\oss 
one has\qss $a_j\qff \not\in\qff A_i$\nnsp,\oss
and\qss
$a_n\qff \not\in\qff A_j$\qss 
for all\qss $j\qff \geqslant\qff k_n$\nnsp.)

(d)\oss \emph{Under the assumptions of\oss \textup{({\halfff}b)},\oss 
show that in order for\qss $m\qff =\qff n$\qss to hold,\oss
it is necessary and sufficient that 
one of the following two cases occurs:}\oss 

({\fff}i\fff)\phantom{i}\oss 
$A_1
\off =\off
\{\qff a_1\fff,\off a_2\fff,\off \ldots\fff,\off a_{n\qff -\qff 1} \qff\}$\nsp,\oss
${A_i
\qff =\qff 
\{\qff a_{i\qff -\qff 1}\fff,\off a_n \qff\}}$\oss
\emph{for}\oss 
$i
\off =\off 
2\fff,\pff \ldots\fff,\pff n$\nsp;\dff

({\fff}ii\fff)\oss 
$n\off =\off k\dff(k\qff -\qff 1)\qff +\qff 1$\nnsp,\oss
\emph{each\dss $A_i$\dss is a set of\qss $k$\qss elements,\oss 
and each element of\qss $E$\qss belongs to exactly $k$ sets $A_i$\nnsp.\oss}

\vspace{6pt}
\myuppar{Remarks.}
Two aspects of this exercise need to be clarified.\oss
First\halfff,\oss
the parts $A_i$ are implicitly assumed to be different from $E$\nnsp.\oss
Second,\oss
the case\qss ({\fff}ii\fff)\qss of the part\dss (d)\qss is expected to hold
only up to renumbering of elements $a_i$ and parts $A_j$\nnsp.\oss

\myuppar{The troubles with the hint\halfff.}
The parts\dss (a)\dss and\dss (b)\dss of this exercise are rather easy,\oss
and there is a hint for the part\dss (c).\oss
But for me this hint turned out to be more of a riddle than of a help.\oss

It would be quite easy to accept and follow the
suggestion to consider the least of the numbers $k_j$\nnsp.\qss
But why it should be $k_n$\nsp?\oss
The phrase\qss 
\emph{``Let\dss $k_n$\dss be the least of the numbers\dss $k_j$\dnsp''}\qss
is fairly hard to interpret\qss
(the expressions used in the French original and 
in the Russian translation have the same meaning\halfff).\qss
The standard usage of\qss \emph{``Let''}\qss
(and of\qss \emph{``Soit''}\qss in French)\qss
in mathematics is to introduce new notations.\qss
But $k_n$ is already defined.\qss

The authors of the solution\qss \cite{q-solution}\qss
found a clever way out\halfff.\oss
They introduce the number $k_n$\dss \emph{before}\pss
introducing other numbers $k_j$\nsp!\oss
This trick helps only partially:\qss
the question\qss \emph{``Why $k_n$\dnsp?''}\qss remains.
The de Bruijn--Erd\"{o}s exposition\qss \cite{db-e}\qss is better\halfff.\qss
They write\qss
\emph{``Assume now that $k_n$ is the smallest $k_i$ \ldots''}.\qff\oss
This\dss is\dss less obscure,\qss 
and amounts to renumbering elements of $E$\nnsp,\oss
but\dss leaves the question\qss \emph{``Why $k_n$\dnsp?''}\oss 
unanswered.\qss

If one manages to put this question aside,\qss
there is another riddle:\qss how the subscripts\qss $i\fff,\pff j$\dnsp,\oss
which are merely marking the points\qss
(and do not even need to be numbers)\qss
may be compared with $k_n$\nnsp,\qss
which is a genuine characteristic of the point marked by $n$\nsp?\oss
Perhaps,\qss this difficulty is encountered only by the 
categorically minded mathematicians;\qss
analysts appear to be quite comfortable with using the values
of a function in its domain of definition.\qss

Here de Bruijn and Erd\"{o}s\qss \cite{db-e}\qss are again doing better\halfff.\oss
They write\qss
\emph{``Assume\oss \ldots\oss that\qss 
$A_1$\nnsp,\qss $A_2$\nnsp,\qss \ldots\dff,\qss $A_{\dff k_n}$\qss
are lines through $a_n$\nnsp''}\qss
(they call the parts $A_i$ lines).\oss
This amounts to renumbering the parts $A_i$\nnsp,\qss 
and one may wonder why renumbering is treated as an assumption.\qss
The trick of the authors of\qss \cite{q-solution}\qss saves the day here for them.\oss
They simply denote the $k_n$\dss lines through $a_n$ by\qss
$A_1$\nnsp,\qss $A_2$\nnsp,\oss \ldots\dff,\qss $A_{\dff k_n}$\qss
and other lines by\qss
$A_{\dff k_n\dff +\dff 1}$\nnsp,\qss $A_{\dff k_n\dff +\dff 2}$\nnsp,\oss
\ldots\dff,\qss $A_{m}$\nnsp.

There is one more riddle in the store.\qss
How one uses the assumption that $k_n$ is the least of the numbers $k_j$
in the proof of the claim in the hint\fff?\oss
One does not\halfff,\qss
this claim is true without it\halfff.

\myuppar{Partially decrypting the hint\halfff.}
Even if one encounters all these troubles and is not aware of the
de Bruijn--Erd\"{o}s paper\qss
({\fff}like me in 1970),\qss
the hint still may be of some help.\oss
The first message is that it is important to know when an element $a_j$
is not in the part $A_i$\nnsp.\oss
Together with the part\qss (b)\qss this suggest that the inequalities\qss
$s_{\fff i}\qff \leqslant\qff k_{\fff j}$\nnsp,\qss
which hold for\qss $a_j\qff \not\in\qff A_i$\dnsp,\oss
should play a key role.\qss

Another message is that the least of the numbers $k_{\fff j}$ should play some role.\qss 
After wasting some time assuming that for a given $u$ 
the number $k_{\fff u}$ is minimal among all numbers $k_{\fff j}$ 
and\dss trying to use this minimality to prove 
something like stated in the hint\halfff,\qss
it is only natural to abandon this assumption and 
consider an arbitrary subscript $u$ such that\qss
$1\qff \leqslant\qff u\qff \leqslant\qff n$\nnsp.\oss

\myuppar{The\qss 1970\qss proof\qss of\qss $m\qff \geqslant\qff n$\nnsp.}
With no more than this limited help from this exercise\qss
(in\dss 1970\dss I definitely understood less than in 2016)\qss
I managed to prove in the early\dss 1970\dss the inequality\qss $m\qff \geqslant\qff n$\nnsp.\oss
Among my schoolmates this qualified as a solution of the Problem\qss M5.\oss
This solution was lost long time ago.\oss
In April of\qss 2016\qss and another time one year later
I attempted to reconstruct this proof\halfff.\oss
In these attempts I encountered the same difficulties as in\dss 1970,\oss
and it is likely that I dealt with them in the same manner\halfff.\oss
At the very least,\qss the resulting proof does not use any tools not known to me at the time,\oss
and does not involve any tricks\qss
(such as the cyclic ordering of some parts $A_i$ by de Bruijn--Erd\"{o}s)\qss
which I was unlikely to discover at the time.\oss
It is presented in Section\qss \ref{solution}\qss below.\oss

The question\qss \emph{``When an equality is possible\fff?''}\qss
was considered by my classmates as too vague to be addressed seriously,\oss 
and this was indirectly admitted by the authors of the solution\qss \cite{q-solution}.\oss
If\qss $m\qff =\qss n$\nnsp,\oss
then\qss (d)\qss easily implies that\qss
$A_i\dff \cap\dff A_j\qff \neq\qff \emptyset$\qss if\qss $i\qff \neq\qff j$\nnsp.\oss
In fact\halfff,\pss proving this property is an almost inevitable part of the proof of\qss (d).\oss
This property means that the set $E$ together with the parts $A_i$ is a\qss
\emph{finite projective plane},\pss
possibly degenerate in the case\qss ({\fff}i\fff)\qss of the part\qss (d).\oss
Therefore,\oss this question amounts to the classification of
finite projective planes and,\qss
to the best of my knowledge,\qss it remains largely open.\oss
See the paper by Ch.\qss Weibel\qss \cite{w}\qss for a survey 
of the state of the art as of\qss 2007,\qss
and\qss \cite{i-planes}\qss for an introduction\qss
(not focusing on the finite case).

\myuppar{\emph{``Kvant''}\qss publishes a solution.}
\emph{``Kvant''}\qss published a solution\qss \cite{q-solution}\qss of
the Problem\qss M5\qss in the August or September of 1970,\oss
close to the beginning of the school year in the USSR\qss
(always September 1).\oss
The editors of the problem section wrote\qss (see\qss \cite{q-solution},\oss p.\qss 49):\qss
\begin{quote}
The letters to editors indicate that this problem is extremely difficult\halfff,\qss
but interesting.\qss
As a matter of fact\halfff,\oss
here we have two problems:\qss 1)\qss prove that\qss $m\qff \geqslant\qff n$\nnsp,\oss
2)\qss when an equality is possible? 

The first problem was completely solved only by\qss \emph{A.\qss Suslin}\qss
from the city of Le\-nin\-grad.\qss
His proof is based on a basic theorem of the linear algebra:\qss
if the number of\dss $n$\dnsp-vectors is greater than $n$\nnsp,\oss
then they are linearly dependent\halfff.\qss

Looking for such a proof will be interesting for whose who are familiar with these notions.\qss
Nobody solved completely the second problem.\qss
Of course,\qss this is not surprising,\qss
since,\qss as it will be explained below,\qss
it can be reduced to a well known,\qss 
but unsolved problem in mathematics.
\end{quote}
Among my schoolmates,\qss these remarks stirred a renewed interest in the problem.\qss
A.\qss Suslin\qss was known as a very strong problem solver and as a
winner of the gold medal at 1967 
International Mathematical Olympiad.\qss
Since only he submitted a complete solution,\qss
the problem had to be really difficult\halfff.\qss
Since he used tools going beyond the school level,\qss
the problem had to be even more difficult.\qss
And this caused a real interest in my unsubmitted 
to the\qss \emph{``Kvant''}\qss solution.\qss
I had an outline as a sparsely filled with formulas sheet of paper\halfff.\qss
One of my schoolmates borrowed this sheet for few days,\qss
and I have not seen it anymore.\qss

But I am not aware of any serious attempt
to study the published solution\qss \cite{q-solution}.\qss
For me it was almost as condensed and obscure as the N.\qss Bourbaki hint\halfff.\oss
The role of the numbering of elements and parts is overemphasized: 
\begin{quote}
Let us pay attention once again to the way we numbered elements and sets.\oss

First of all,\qss $k_{\fff n}$\dss is the least of the numbers\qss
$k_{\fff 1}\fff,\pff k_{\fff 2}\fff,\pff \ldots\fff,\pff k_{\fff n\dff -\dff 1}$\qss
({\fff}sic\fff!\oss --\qss \emph{N.I.}\trf).\oss \ldots
\end{quote}
See\qss \cite{q-solution},\oss p.\qss 51.\oss
And I always disliked random numerical examples,\qss
which are supposed to help the reader and are
extensively used in\qss \cite{q-solution}.\qss
I must admit that I did not even look at the last two pages of\qss \cite{q-solution}\qss
before writing these comments,\pss
and,\pss in particular\halfff,\pss before writing down the proof\dss in the next section.\oss
Surprisingly,\oss it turned out that the proof\qss \cite{q-solution}\qss contains a gap:\oss
it is mentioned that\qss $k_{\fff n}\qff =\qff 2$\qss in the situation
described in the case\qss ({\fff}i\fff)\qss of the Bourbaki exercise,\oss
but no proof that this is the only possibility is even attempted.\oss

\mysection{A\qss solution\qss of\qss the\qss N.\qss Bourbaki\qss exercise}{solution}

\vspace*{6pt}
\myuppar{The terminology and notations.}
In contrast with\dss N.\qss Bourbaki\dss and with the\qss \emph{``Kvant''},\pss
I\dss have no reasons to hide the geometric content of this result\halfff.\qss
Following de Bruijn and\dss Erd\"{o}s,\oss I will call the elements of $E$\dss \emph{points}\qss
and\dss the sets $A_i$\dss \emph{lines}.\pss
Since the lines are assumed to be proper subsets of $E$\nnsp,\pss
every point is contained in at least $2$ lines.\oss
Indeed,\oss if a point is contained in only one line,\pss
then all points are contained in this line,\oss
i.e.\qss it is not a proper subset\halfff.\oss

It is convenient 
to explicitly introduce a counterpart to the set $E$ of points,\oss
namely the set of lines\oss
$\mathcal{L}
\off =\off
\{\dff A_1\fff,\pff A_2\fff,\pff \ldots\fff,\pff A_m \dff\}$\nnsp.\oss
If\dss the case\qss ({\fff}i\fff)\qss of the part\qss (d)\qss of the 
Bourbaki exercise occurs,\pss
up to renumbering of points and lines,\pss
then the pair\qss $(\dff E\fff,\pff \mathcal{L}\dff)$\qss 
is called a\qss \emph{near-pencil}.\oss
If the case\qss ({\fff}ii\fff)\qss of the part\qss (d)\qss occurs,\oss
then\qss $(\dff E\fff,\pff \mathcal{L}\dff)$\qss is called a\qss \emph{projective plane}.\oss

I also do not see any reason to follow the outdated fashion of using\qss
\emph{numerical}\pss indices\qss (i.e.\qss subscripts),\qss
which amounts to ordering objects even when their order is irrelevant\halfff.\qss
Instead of this,\qss
for every point $z$ 
we will denote by $k_{\fff z}$ the number of lines containing $z$\nnsp,\pss
and for every line $l$ 
we will denote by $s_{\fff l}$ the number of points in $l$\dnsp,\oss
i.e. the number of elements of the set $l$\dnsp.

\myuppar{The part\qss (a)\qss of the Bourbaki exercise.}
With the above notations the part\qss (a)\qss takes the form\vspace*{2pt}
\begin{equation}
\label{sums}
\quad
\sum_{l\dff \in\dff \mathcal{L}}\qff s_{\fff l}
\off\off =\off\off
\sum_{\qff z\dff \in\dff E}\qff k_{\fff z}\dff.
\end{equation}

\vspace*{-10pt}
after interchanging the sides.\oss
This immediately follows from counting in two different ways the pairs\dss
$(z\fff,\pff l\dff)\qff \in\qff E\times \mathcal{L}$\qss such that\dss $z\qff \in\qff l$\nnsp.\oss

\myuppar{The part\qss ({\halfff}b\halfff)\qss of the Bourbaki exercise.}
For the rest of the paper we will assume that 
the assumption of the part\qss 
({\halfff}b\halfff)\qss 
holds,\oss
i.e.\qss that for every pair of distinct points there is exactly one line containing both of them.\oss
If a line $l$ contains\qss $\leqslant\qff 1$\qss points,\oss 
then removing $l$ from the set of lines
does not affects this assumption,\oss 
and at the same time decreases number of lines by $1$\nnsp.\oss
Hence we may assume for the rest of the paper that
every line contains at least $2$ points.\oss

With the above notations
the part\qss 
({\halfff}b\halfff)\qss 
takes the form\vspace*{2pt}
\begin{equation*}
\quad
\mbox{ If }\quad
z\qff \not\in\qff l\fff,\quad
\mbox{ then }\quad
s_{\fff l}\qff \leqslant\qff k_{\fff z}\qff.
\end{equation*}

\vspace*{-10pt}
We will call these inequalities
the\qss \emph{\dbe {\dff}inequalities}.\oss

In order to prove the \dbe inequalities,\oss
suppose that\qss  $z\qff \not\in\qff l$\nnsp.\oss
Then for every\dss $y\qff \in\qff l$\dss there is a unique line
containing\dss $\{\dff z\fff,\pff y \dff\}$\dss and it is different from $l$
because\qss  $z\qff \not\in\qff l$\dnsp.\oss
These lines are pairwise distinct because if\qss
$y\fff,\pff y'\qff \in\qff l$\qss and\qss $y\qff \neq\qff y'$\dnsp,\oss
then $l$ is the only line containing\dss $\{\dff y\fff,\pff y' \dff\}$\nnsp.\oss
There are\dss is $s_{\fff l}$ such lines and all of them contain $z$\nnsp;\oss
therefore\qss $s_{\fff l}\qff \leqslant\qff k_{\fff z}$\nnsp.\oss

\myuppar{Lines through an arbitrary point\halfff.}
Let\qss $u\qff \in\qff E$\qss be an arbitrary point\halfff,\oss
let\qss $p\qff =\qff k_{\fff u}$\qss be the number of lines containing $u$\nnsp,\oss
and let $\mathcal{U}$ be the set of these lines.\oss

By the definition of $\mathcal{U}$\dnsp,\oss
if a line $l$\dss is not in\dss $\mathcal{U}$\dnsp,\oss
then\qss $u\qff \not\in\qff l$\dnsp.\oss
For every\qss 
$l\qff \not\in\qff \mathcal{U}$\qss
we have the \dbe inequality\qss
$s_{\fff l}\qff \leqslant\qff k_{\fff u}$\nnsp.\oss
By summing all these inequalities and taking into account that there
are\qss $m\qff -\qff p$\qss lines not\dss belonging\dss to $\mathcal{U}$\dnsp,\oss
we see that\vspace*{2pt}
\begin{equation}
\label{s-upper-estimate}
\quad
\sum_{l\qff \not\in\qff \mathcal{U}} s_{\fff l}
\off\off \leqslant\off\off 
(m\qff -\qff p)\dff k_u\qff.
\end{equation}

\vspace*{-10pt} 
Since every set of the form\dss 
$\{\dff u\fff,\pff y \trf\}$\dss with\dss
$y\qff \neq\qff u$\dss 
is contained in one and only one line,\oss
the sets\qss
$l\dff \smallsetminus\dff \{\dff u\trf\}$\qss
with\qss $l\qff \in\qff \mathcal{U}$\qss
are pairwise disjoint and form a partition of\qss
$E\dff \smallsetminus\dff \{\dff u\trf\}$\nnsp.\oss
Since we assumed that\qss $s_l\qff \geqslant\qff 2$\qss for all lines $l$\dnsp,\oss
all these sets are non-empty.\oss
Let $U$ be a set of representatives of these sets.\oss
In other terms,\qss $U$ is contained\dss in\qss
$E\dff \smallsetminus\dff \{\dff u\trf\}$\qss
and intersects every set\qss
$l\dff \smallsetminus\dff \{\dff u\trf\}$\qss
with\qss $l\qff \in\qff \mathcal{U}$\qss
in exactly $1$ point\halfff.\oss
In particular\halfff,\oss $U$ consists of exactly $p$ points.\oss

If\qss 
$(\fff l\fff,\pff z\fff)\qff \in\qff \mathcal{U}\dff \times\dff U$\qss 
and\qss 
$z\qff \not\in\qff l$\nnsp,\oss
then the \dbe inequality\qss
$s_{\fff l}\qff \leqslant\qff k_{\fff z}$\qss holds.\oss
There are $p\fff(p\qff -\qff 1)$ of such pairs\qss
$(\fff l\fff,\pff z\fff)$\qss and hence\qss
$p\fff(p\qff -\qff 1)$\qss of such inequalities.\oss
For each\dss $l\qff \in\qff \mathcal{U}$\dss
the number $s_{\fff l}$ occurs\dss $p\qff -\qff 1$\dss times in the left hand sides
of them,\pss
and for each\dss $z\qff \in\qff U$\dss the number $k_{\fff z}$
occurs\dss $p\qff -\qff 1$\dss times in the right hand sides.\oss
Hence the sum of all these inequalities is\vspace*{2pt}
\begin{equation}
\label{sum-all-pairs}
\quad
\sum_{l\qff \in\qff \mathcal{U}} (p\qff -\qff 1)\dff s_{\fff l}
\off\off \leqslant\off\off 
\sum_{z\qff \in\qff U} (p\qff -\qff 1)\dff k_{\fff z}\qff.
\end{equation}

\vspace*{-10pt}
After dividing\pss (\ref{sum-all-pairs})\pss 
by\pss $p\qff -\qff 1$\pss we get\vspace*{2pt}
\begin{equation}
\label{sum-divided}
\quad
\sum_{l\qff \in\qff \mathcal{U}} s_{\fff l}
\off\off \leqslant\off\off 
\sum_{z\qff \in\qff U} k_{\fff z}\qff.
\end{equation}

\vspace*{-10pt}
Now it is only natural to take the sum of the inequalities\qss 
(\ref{s-upper-estimate})\qss and\qss (\ref{sum-divided})\qss 
and conclude that\vspace*{2pt}
\begin{equation}
\label{sum-all-s-estimate}
\quad 
\sum_{l\qff \in\qff \mathcal{L}} s_{\fff l}
\off\off \leqslant\off\off
(m\qff -\qff p)\dff k_{\fff u}
\off +\off 
\sum_{z\qff \in\qff U} k_{\fff z}\qff.
\end{equation}

\vspace*{-10pt} 
The left hand side of the inequality\qss (\ref{sum-all-s-estimate})\qss 
is the same as the left hand side of the equality\qss (\ref{sums}).\qss
The right hand side of\qss (\ref{sum-all-s-estimate})\qss 
can be compared with the right hand side of the equality\qss (\ref{sums})\pss 
if\qss $k_{\fff u}$\qss is\dss the least among the numbers\qss $k_{\fff z}$\qss
and\qss $m\qff \leqslant\qff n$\nnsp.\oss

\myuppar{Proof\dss of\dss the\dss inequality\qss $m\qff \geqslant\qff n$\nnsp.}
Now we are ready to do the part\qss (c)\qss of the Bourbaki exercise.\oss
Let\qss $u\qff \in\qff E$\qss be a point such that $k_{\fff u}$ 
is the least of the numbers $k_{\fff z}$
over all points\qss $z\qff \in\qff E$\nnsp.\oss
Then\vspace*{2pt}
\begin{equation}
\label{k-lower-estimate}
\quad
(m\qff -\qff p)\dff k_{\fff u}
\off\off \leqslant\off\off
\sum_{z\qff \in\qff Y} k_{\fff z}\qff.
\end{equation}

\vspace*{-10pt}
for every subset\qss $Y\qff \subset\qff E$\qss consisting of\qss $m\qff -\qff p$\qss points.\oss

Suppose that\qss $m\qff \leqslant\qff n$\nnsp.\oss
Then the subset $Y$ can be chosen to be disjoint from $U$\qss
({\fff}because $U$ consists of $p$ points).\oss
Let us choose an arbitrary $Y$ disjoint from $U$ and\dss let\qss
$Z\off =\off Y\qff \cup\qff U$\dnsp.\oss
Then $Z$ is a subset of $E$ consisting of $m$ points and
the inequalities\qss 
(\ref{sum-all-s-estimate})\qss and\qss (\ref{k-lower-estimate})\qss 
imply that\vspace*{2pt}
\begin{equation}
\label{final-estimate}
\quad
\sum_{l\qff \in\qff \mathcal{L}} s_{\fff l}
\off\off \leqslant\off\off
\sum_{z\qff \in\qff Z} k_{\fff z}
\off\off \leqslant\off\off
\sum_{z\qff \in\qff E} k_{\fff z}\qff,
\end{equation}

\vspace*{-10pt}
where the last inequality is strict unless\qss $Z\off =\off E$\nnsp.\oss
In view of\qss (\ref{sums})\qss this inequality cannot be strict
and hence\qss $Z\off =\off E$\qss and\qss $m\qff =\qff n$\nnsp.\oss
Since\qss $m\qff \leqslant\qff n$\qss implies that\qss $m\qff =\qff n$\nnsp,\oss
we see that\qss $m\qff \geqslant\qff n$\nnsp.\oss

\myuppar{The case\qss $m\qff =\qff n$\nnsp.}
After the work done in the proofs of\qss (a),\qss ({\halfff}b\halfff),\qss and\qss (c),\qss
the part\qss (d\halfff)\qss nearly proves itself\halfff.\qss
As we will see,\qss
in this case all inequalities\qss 
(\ref{s-upper-estimate})\qss --\qss (\ref{final-estimate})\qss
are,\oss in fact\halfff,\oss equalities.

By\qss (\ref{sums})\qss the leftmost and the rightmost sums in\qss
(\ref{final-estimate})\qss are equal.\oss
It follows that\qss $Z\off =\off E$\qss and hence\qss
$Y\off =\off E\dff \smallsetminus\dff U$\dnsp.\oss
Moreover\halfff,\oss 
the sides of each of the inequalities\qss 
(\ref{sum-all-s-estimate})\qss and\qss (\ref{k-lower-estimate})\qss
are equal.\oss 
Since $k_{\fff u}$ is the least of the numbers $k_{\fff z}$\nnsp,\oss
the equality of the sides of\qss (\ref{k-lower-estimate})\qss
implies that\vspace*{2pt}
\begin{equation}
\label{ku-k-not-u}
\quad
k_{\fff u}
\off =\off
k_{\fff z}
\quad\mbox{ for\qss all }\quad
z\qff \in\qff Y\off =\off E\dff \smallsetminus\dff U\qff.
\end{equation}

\vspace*{-10pt}
The fact that the sides of\qss
(\ref{sum-all-s-estimate})\qss
are equal implies that the sides of each of the inequalities\qss
(\ref{s-upper-estimate})\qss and\qss (\ref{sum-divided})\qss
are equal also.\oss
The equality of the sides of\qss (\ref{s-upper-estimate})\qss 
implies that\qss\vspace*{2pt}
\begin{equation}
\label{s-not-ku}
\quad
s_{\fff l}
\off =\off 
k_{\fff u}
\quad\mbox{ for\qss all }\quad
l\qff \not\in\qff \mathcal{U}\qff. 
\end{equation}

\vspace*{-10pt}
Since the sides of\qss (\ref{sum-divided})\qss are equal,\oss
the sides of\qss (\ref{sum-all-pairs})\qss
are also equal.\oss
Since\qss (\ref{sum-all-pairs})\qss is the sum of the inequalities\qss
$s_{\fff l}\qff \leqslant\qff k_{\fff z}$\qss 
over all pairs\qss 
$(\fff l\fff,\pff z\fff)\qff \in\qff \mathcal{U}\dff \times\dff U$\qss 
such that\qss $z\qff \not\in\qff l$\nnsp,\oss
the equality of the sides of\qss (\ref{sum-all-pairs})\qss
implies that\qss
$s_{\fff l}\qff =\qff k_{\fff z}$\qss
for all\qss 
$(\fff l\fff,\pff z\fff)\qff \in\qff \mathcal{U}\dff \times\dff U$\qss 
such that\qss $z\qff \not\in\qff l$\nnsp.\oss
Equivalently,\vspace*{2pt}
\begin{equation}
\label{s-k-uu}
\quad
s_{\fff l}
\off =\off 
k_{\fff z}
\quad\mbox{ if }\quad
l\qff \in\qff \mathcal{U} 
\quad\mbox{ and }\quad
z\qff \in\qff U\qff \smallsetminus\qff l\qff.
\end{equation}

\vspace*{-10pt}
The rest of the proof splits into two subcases depending on\dss if\qss 
$p\qff =\qff 2$\qss or\qss $p\qff \geqslant\qff 3$\nnsp.

\myuppar{The subcase\qss $p\qff =\qff 2$\nnsp.}
In this case\qss $\mathcal{U}\off =\off \{\trf l\fff,\pff l' \qff\}$\qss
for some\qss $l\fff,\pff l'$\qss and\dss hence\qss
$E\qff =\qff l\qff \cup\qff l'$\dnsp.\oss
It follows that every line different from\qss $l\fff,\pff l'$\qss
contains only $2$ points,\oss namely
the points of its intersection with the lines $l\fff,\pff l'$\dnsp.\qff\oss
If\oss $s_{\fff l}\fff,\off s_{\fff l'}\pff \geqslant\pff 3$\nnsp,\oss
then there are at least $4$ points\qss $z\qff \neq\qff u$ and the part\qss 
({\halfff}b)\qss implies that\qss 
$k_{\fff z}\qff \geqslant\qff 3$\qss for every\qss $z\qff \neq\qff u$\nnsp.\oss
On the other hand,\qss 
(\ref{ku-k-not-u})\qss implies that\vspace*{2pt}
\[
\quad
k_{\fff z}
\off =\off 
k_{\fff u}
\off =\off 
p
\off =\off 
2
\]

\vspace*{-10pt}
for every\qss $z\qff \not\in\qss U$\nnsp.\oss
But $U$ consists of only two points and hence\qss
$k_{\fff z}\qff \geqslant\qff 3$\qss for no more than two points $z$\nnsp.\oss
The contradiction shows that
either\qss $s_{\fff l}\qff =\qff 2$\qss or\qss $s_{\fff l'}\qff =\qff 2$\nnsp.\oss
We may assume that\qss $s_{\fff l}\qff =\qff 2$\nnsp.\oss
Then\qss $l\qff =\qff \{\dff u\fff,\pff a \trf\}$\qss for some\qss $a\qff \in\qff E$\qss
and every line different from\qss $l\fff,\pff l'$\qss
has the form\qss $\{\dff a\fff,\pff z \dff\}$\qss with\qss
$z\qff \in\qff l'\dff \smallsetminus\dff \{\dff u\trf\}$\nnsp.\oss
It\dss follows that\qss $(\fff E\fff,\pff \mathcal{L}\fff)$\qss 
is a near-pencil.\oss

\myuppar{The subcase\qss $p\qff \geqslant\qff 3$\nnsp.}
The set $U$ is a set of representatives of the sets\qss
$l\dff \smallsetminus\dff \{\fff u\trf\}$\qss
with\qss $l\qff \in\qff \mathcal{U}$\dnsp.\oss
For any two lines\qss $l\fff,\pff l'\qff \in\qff \mathcal{U}$\qss
the assumption\qss $p\qff \geqslant\qff 3$\qss
implies that there exists a point\qss $z\qff \in\qff U$\qss
such that\qss $z\qff \not\in\qff l\fff,\pff l'$\dnsp.\qff\oss
If\qss $z$\qss is such a point\halfff,\oss
then\qss (\ref{s-k-uu})\qss implies that
\[
\quad
s_{\fff l}\off =\off k_{\fff z}\off =\off s_{\fff l'}\qff.
\]
Similarly,\oss if\qss $z\fff,\pff z'\qff \in\qff U$\dnsp,\oss
then there exists a line\qss $l\qff \in\qff \mathcal{U}$\qss
such that\qss $z\fff,\pff z'\qff \not\in\qff l$\qss and hence
\[
\quad
k_{\fff z}\off =\off s_{\fff l}\off =\off k_{\fff z'}\qff.
\]
It follows that in the subcase\qss $p\qff \geqslant\qff 3$\qss
all numbers\qss $s_{\fff l}$\nnsp,\qss $k_{\fff z}$\qss
with\qss $l\qff \in\qff \mathcal{U}$\qss
and\qss $z\qff \in\qff U$\qss are equal.\oss
Since\qss $k_{\fff u}\qff =\qff p\qff \geqslant\qff 3$\qss
is the smallest of the numbers $k_{\fff z}$ over all\qss $z\qff \in\qff E$\nnsp,\oss
it follows that\qss
\[
\quad
s_{\fff l}
\off =\off
k_{\fff z}
\off \geqslant\off 3
\quad\mbox{ for\qss all }\quad
l\qff \in\qff \mathcal{U}\fff,\quad
z\qff \in\qff U.
\]
Let\qss $l\qff \in\qff \mathcal{U}$\dnsp,\oss
and let $y$ be the unique element of $U$ contained in $l$\dnsp.\oss 
Since\qss $s_{\fff l}\qff \geqslant\qff 3$\nnsp,\oss 
there exists a point\qss $x\qff \in\qff l$\qss
not equal to\qss $u\fff,\pff y$\nnsp.\oss
We can replace in $U$ the point $y$ by the point $x$
and get a new set of representatives $U'$\dnsp.\oss
Then all previous results apply to $U'$ in the role of $U$\dnsp.\oss
In particular\halfff,\oss
$k_{\fff x}
\off =\off 
k_{\fff z}$\oss
for\qss all\oss 
$z\qff \in\qff U\dff \smallsetminus\dff l
\off =\off
U'\dff \smallsetminus\dff l$\oss
and hence
\[
\quad
k_{\fff x}
\off =\off 
k_{\fff z}
\quad\mbox{ for\qss all }\quad
z\qff \in\qff U.
\]
On the other hand,\oss $x\qff \not\in\qff U$\qss and hence\qss 
$k_{\fff x}\qff =\qff k_{\fff u}$\qss by\qss
(\ref{ku-k-not-u})\qss applied to the original set $U$\dnsp.\oss
At the same time\qss (\ref{ku-k-not-u})\qss
implies that\qss
$k_{\fff u}\qff =\qff k_{\fff z}$\qss
for all\qss $z\qff \not\in\qff U$\qss
and hence\qss
\[
\quad
k_{\fff x}
\off =\off 
k_{\fff z}
\quad\mbox{ for\qss all }\quad
z\qff \not\in\qff U.
\] 
It\dss follows that all numbers $k_{\fff z}$ are equal.\oss
At the same time by\qss (\ref{s-not-ku})\qss and\qss (\ref{s-k-uu})\qss
every $s_{\fff l}$\dss is equal to some $k_{\fff z}$\nnsp.\oss
It follows that all numbers\qss $s_{\fff l}$\nnsp,\qss $k_{\fff z}$\qss
with\qss $l\qff \in\qff \mathcal{L}$\qss
and\qss $z\qff \in\qff E$\qss  are equal.\oss
It remains to apply the following lemma.\oss

\myuppar{Lemma\qss 1.}
\emph{If\trs all\dss the numbers\qss $s_{\fff l}$\nnsp,\qss $k_{\fff z}$\qss
are equal,\oss
then\qss $(\dff E\fff,\pff \mathcal{L}\dff)$\qss is a projective plane.}

\myuppar{Proof\halfff.}
Let $k$ be the common value of the numbers\qss $s_{\fff l}$\nnsp,\qss $k_{\fff z}$\nnsp,\oss
and\dss let\qss $y\qff \in\qff E$\nnsp.\oss
The sets\qss
$l\dff \smallsetminus\dff \{\dff y\trf\}$\qss
with\qss $y\qff \in\qff l$\qss
are pairwise disjoint and form a partition of\qss
$E\dff \smallsetminus\dff \{\dff y\trf\}$\nnsp.\oss
Each of them consists of\qss 
\[
\quad
s_{\fff l}\qff -\qff 1\off =\off k\qff -\qff 1
\]
points,\oss and there are\qss $k_{\fff y}\qff =\qff k$\qss such sets.\oss
It follows that the number $n$ of elements of $E$ is equal to\qss
$k\fff(k\qff -\qff 1)\qff +\qff 1$\nnsp.\oss
Therefore,\oss $(\dff E\fff,\pff \mathcal{L}\dff)$\qss 
is a projective plane.\oss  \eproof

\myuppar{Remarks.}
A key step of this solution and the solution\qss \cite{q-solution}\qss
differ from the de-Bruijn--Erd\"{o}s paper in the same way:\qss
the cyclic order argument of \dbe 
(see Section\qss \ref{dbe-proof})\qss
is replaced by the inequalities\qss (\ref{sum-all-pairs})\qss and\qss (\ref{sum-divided}).\oss

\mysection{The\qss de Bruijn--Erd\"{o}s\qss proof}{dbe-proof}

\vspace*{12pt}
Since the proof presented in Section\qss \ref{solution}\qss grow out of a
summary of the de Bruijn--Erd\"{o}s proof\halfff,\oss
albeit not quite understood,\oss
it is not surprising that the two proofs have a lot in common.\oss
In the following exposition of the de Bruijn--Erd\"{o}s proof 
we will use the notations of Section\qss \ref{solution}\qss and 
will refer to Section\qss \ref{solution}\qss for the arguments 
which differ from\qss \cite{db-e}\qss only in the notations and the amount of details.\oss
The Bruijn--Erd\"{o}s paper is concise on the border of being cryptic.\oss

The de Bruijn--Erd\"{o}s proof begins with the parts\qss (a)\qss and\qss ({\halfff}b\halfff)\qss
of the Bourbaki exercise.\oss
After this \dbe introduce $k_{\fff u}$ as the smallest among all numbers $k_{\fff z}$\qss
(and denote it by $k_{\fff n}$\nnsp).\oss
Then \dbe observe that 
it can assumed that every line contains at least two points.\oss
Following the notations of Section\qss \ref{solution},\oss
let us denote by $\mathcal{U}$ the set of all lines containing $u$\nnsp.\oss
By the \dbe inequalities\qss
$s_{\fff l}\qff \leqslant\qff k_{\fff u}$\qss 
for every\qss 
$l\qff \not\in\qff \mathcal{U}$\dnsp.\oss
The inequalities\qss (\ref{s-upper-estimate})\qss
and\qss (\ref{k-lower-estimate})\qss follows.\oss
The following argument plays a role similar to the role of the inequality\qss
(\ref{sum-all-pairs}).

\myuppar{The cyclic order argument\halfff.}
Let\oss
$l_{\fff 1}\fff,\pff l_{\fff 2}\fff,\pff \ldots\fff,\pff l_{\fff p}$\oss
be a cyclically ordered list of elements of $\mathcal{U}$\dnsp.\oss
We treat the subscripts\qss $1\fff,\pff 2\fff,\pff \ldots\fff,\pff p$\qss
as integers\dnsp\ $\mod\dss p$\dnsp.\oss
For each\oss 
$i
\off =\off
1\fff,\pff 2\fff,\pff \ldots\fff,\pff p$\oss
let us choose some point\qss $a_{\fff j}\qff \in\qff l_{\fff j}\dff \smallsetminus\dff \{\dff u \dff\}$\qss
and\dss let $U$ be the set of these points.\oss
Let 
\[
\quad
s_{\fff i}
\off =\off 
s_{\fff l_{\fff i}}
\hspace*{1.5em}\mbox{ and }\hspace*{1.5em}
k_j
\off =\off 
k_{\fff a_{\fff j}}\qff.
\]
Since\qss $a_{\fff i\dff +\dff 1}\qff \not\in\qff l_{\fff i}$\nnsp,\oss 
by the \dbe inequalities
$s_{\fff i}\qff \leqslant\qff k_{\fff i\dff +\dff 1}$\qss
for all\qss 
$i\off =\off
1\fff,\pff 2\fff,\pff \ldots\fff,\pff p$\dnsp,\oss
i.e.\vspace*{2pt}
\begin{equation}
\label{s-k-cycle}
\quad
s_{\fff 1}\off \leqslant\off k_{\fff 2}\dff,\quad
s_{\fff 2}\off \leqslant\off k_{\fff 3}\dff,\quad
\ldots\dff,\quad
s_{\fff p}\off \leqslant\off k_{\fff 1}\dff.
\end{equation}

\vspace*{-10pt}
By summing the inequalities\qss (\ref{s-k-cycle})\qss one concludes that\oss
\begin{equation}
\label{sum-to-p}
\quad
\sum_{j\qff =\qff 1}^p s_j
\off \leqslant\off
\sum_{j\qff =\qff 1}^p k_j\qff.
\end{equation}
The inequality\qss (\ref{sum-to-p})\qss is nothing else but 
another form of\qss 
(\ref{sum-divided}).\oss 
The arguments of Section\qss \ref{solution}\qss
show that\qss (\ref{sum-to-p})\qss 
implies that\qss
$m\qff \geqslant\qff n$\nnsp.\oss 
In fact\halfff,\oss de Bruijn and Erd\"{o}s do not bother to
write down even the inequality\qss (\ref{sum-to-p}),\oss
to say nothing about other details presented in Section\qss \ref{solution}.

\myuppar{The case\qss $m\qff =\qff n$\nnsp.}
In view of the equality\qss (\ref{sums})\qss
in this case the left hand and the right hand sides of 
the inequality\qss (\ref{sum-to-p})\qss 
are equal.\oss
Together with\qss (\ref{s-k-cycle})\qss this implies that
\begin{equation}
\label{cycle-of-equalities}
\quad
s_{\fff 1}\off =\off k_{\fff 2}\dff,\quad
s_{\fff 2}\off =\off k_{\fff 3}\dff,\quad
\ldots\dff,\quad
s_{\fff p}\off =\off k_{\fff 1}\qff.
\end{equation}
Similarly,\oss in this case
the left hand and the right hand sides of the inequality\qss
(\ref{k-lower-estimate})\qss 
are equal.\oss
Since\qss $m\qff =\qff n$\nnsp,\oss
one can take\qss
$Y\qff =\qff E\dff \smallsetminus\dff U$\qss in\qss (\ref{k-lower-estimate}).\oss
It follows that\qss
$k_{\fff u}
\qff =\qff
k_{\fff z}$\qss
for all\qss
$z\qff \in\qff E\dff \smallsetminus\dff U$\nnsp.\oss
Finally,\oss the left hand and the right hand sides of the inequality\qss
(\ref{s-upper-estimate})\qss
are equal.\oss
It follows that\qss
$s_{\fff l}
\qff =\qff 
k_{\fff u}$\qss
for all\qss
$l\qff \not\in\qff \mathcal{U}$\dnsp.\oss 
By combining the last two observations,\oss
we see that\qss\vspace*{2pt}
\begin{equation*}
\quad
s_{\fff l}
\off =\off 
k_{\fff z}
\quad\mbox{ for\qss all }\quad
l\qff \in\qff \mathcal{L}\dff \smallsetminus\dff \mathcal{U}\fff,\quad
z\qff \in\qff E\dff \smallsetminus\dff U\dff.
\end{equation*}

\vspace*{-10pt}
Since\qss $m\qff =\qff n$\nnsp,\oss both sets\qss
$\mathcal{L}\dff \smallsetminus\dff \mathcal{U}$\qss
and\qss
$E\dff \smallsetminus\dff U$\qss
consist of\qss $n\qff -\qff p$\qss elements.\oss
It follows that one can number the points and lines in
such a way that\qss (in the notation of the Bourbaki exercise)
\[
\quad
s_{\fff 1}\off =\off k_{\fff 1}\fff,\quad
s_{\fff 2}\off =\off k_{\fff 2}\fff,\quad
\ldots\fff,\quad
s_{\fff n}\off =\off k_{\fff n}\qff.
\]
As the next step,\oss let us renumber the points 
and lines once more and assume that
\begin{equation}
\label{k-order}
\quad
k_{\fff 1}\off \geqslant\off
k_{\fff 2}\off \geqslant\off
\ldots\off
\geqslant\off k_{\fff n}\qff.
\end{equation}
The rest of the proof splits into two subcases depending on\dss if\dss 
$k_{\fff 1}\qff >\qff k_{\fff 2}$\dss or not\halfff.

\myuppar{The subcase\qss $k_{\fff 1}\qff >\qff k_{\fff 2}$\nnsp.}
In this case\qss
$s_{\fff 1}\qff =\qff k_{\fff 1}\qff >\qff k_{\fff i}$\qss
for all\qss $i\qff \geqslant\qff 2$\nnsp.\oss
By the \dbe inequalities this implies that\qss 
$a_{\fff i}\qff \in\qff A_{\fff 1}$\qss for all\qss $i\qff \geqslant\qff 2$\nnsp.\oss
It\dss follows that\qss $(\dff E\fff,\pff \mathcal{L}\dff)$\qss is a near-pencil.\oss

\myuppar{The subcase\pss $k_{\fff 1}\off =\off k_{\fff 2}$\nnsp.}
Suppose that\qss $k_{\fff j}\qff <\qff k_{\fff 1}\qff =\qff k_{\fff 2}$\qss
for some $j$\nnsp.\oss
By the \dbe inequalities $a_{\fff j}$ belongs to the both lines
$A_{\fff 1}$ and $A_{\fff 2}$\nnsp.\oss
This is possible for only one point\halfff,\oss
namely the point of the intersection of the lines
$A_{\fff 1}$ and $A_{\fff 2}$\nnsp.\oss
In view of\qss (\ref{k-order}),\oss
this may happen only if\oss
\[
\quad
k_{\fff 1}
\off =\off
k_{\fff 2}
\off =\off
\ldots
\off =\off 
k_{\fff n\dff -\dff 1}
\off >\off
k_{\fff n}
\]
and hence\oss
$
s_{\fff j}
\off =\off
k_{\fff j}
\off >\off
k_{\fff n}
\off \geqslant\off
2$\oss
for all\qss $j\qff \neq\qff n$\nnsp.\oss
It follows that\oss
$s_{\fff j}\off \geqslant\off 3$\oss
if\qss $j\qff <\qff n$\nnsp.\oss
In particular\halfff,\oss all $k_{\fff n}$ lines containing $a_{\fff n}$
consist of\qss $\geqslant\qff 2$\qss points and all except\halfff,\pss perhaps,\pss
the line $A_{\fff n}$\nnsp,\oss consist of\qss $\geqslant\qff 3$\qss points.\oss
Therefore one can choose $2$ points\qss $x\fff,\pff y\qff \neq\qff a_{\fff n}$\qss
on one of these lines,\oss
and a point\qss $z\qff \neq\qff a_{\fff n}$\qss on some other line.\oss
Let\qss $l_{j}\fff,\pff l_{j'}$\qss be the lines containing the pairs\qss
$\{\dff x\fff,\fff z \dff\}$\qss and\qss
$\{\dff y\fff,\fff z \dff\}$\qss respectively.\oss
Then\qss $j\qff \neq\qff j'$\pss and\pss 
$a_{\fff n}\off \not\in\off l_j\fff,\pff l_{j'}$\nnsp.\oss
Hence the \dbe inequalities imply that\qss
$s_{\fff j}\fff,\pff s_{\fff j'}\off \leqslant\off k_{\fff n}$\nnsp,\oss
contrary to the fact that\qss 
$s_{\fff j}\off >\off k_{\fff n}$\oss if\oss $j\off \neq\off n$\nnsp.\oss
The contradiction shows that all numbers $k_{\fff j}$ are equal,\oss
and hence all numbers\qss $s_{\fff i}\fff,\pff k_{\fff j}$\qss
are equal.\oss
Now the observation at the end of Section\qss \ref{solution}\qss
implies that\qss $(\dff E\fff,\pff \mathcal{L}\dff)$\qss is a projective plane.

\myuppar{Intersection of lines.}
After the proof is completed,\oss \dbe point out
that in the subcase\qss $k_{\fff 1}\qff =\qff k_{\fff 2}$\qss
of the case\qss $m\qff =\qff n$\qss
every two lines intersect\halfff.\oss
Indeed,\oss if\qss $l'\fff,\pff l''$\qss are two disjoint lines
and\qss $a\qff \in\qff l''$\dnsp,\oss
then there are $s_{\fff l'}$ lines
containing $a$ and intersecting $l'$\dnsp,\oss
and still one more line,\oss namely $l''$\dnsp,\oss containing $a$\nnsp.\oss
Therefore\qss $k_{\fff a}\qff \geqslant\qff s_{\fff l'}\qff +\qff 1$\nnsp,\oss
contrary to the fact all numbers\qss $k_{\fff z}\fff,\pff s_{\fff l}$\qss are equal.\oss
In fact\halfff,\oss every two lines obviously intersect 
in the subcase\qss $k_{\fff 1}\qff >\qff k_{\fff 2}$\qss also.

\myuppar{Why $k_{\fff n}$\nsp?}
Now it is clear why the smallest of the numbers $k_{\fff z}$ is denoted by $k_{\fff n}$\nnsp.\oss
The number $k_{\fff n}$ is indeed the smallest 
if the points are ordered in such a way that\qss (\ref{k-order})\qss holds.\oss
At the same time\qss (\ref{k-order})\qss 
plays almost no role in the proof\halfff.\oss
One may speculate that\qss (\ref{k-order})\qss 
and notation $k_{\fff n}$ for the smallest of the numbers $k_{\fff z}$
are remnants of an earlier approach to the theorem.\oss

\mysection{From\qss de Bruijn--Erd\"{o}s\qss to\qss
systems\qss of\qss distinct\qss representatives}{reps}

\vspace*{6pt}
\myuppar{The cyclic order argument and systems of distinct representatives.}
The key step of the de Bruijn--Erd\"{o}s proof is the cyclic order argument
used to prove the inequality\qss (\ref{sum-to-p})\qss 
and the equalities\qss (\ref{cycle-of-equalities})\qss in the case\qss
$m\qff =\qff n$\nnsp.\oss
Ultimately,\oss 
the cyclic order argument is based on the fact that\qss
$a_{\fff i\dff +\dff 1}\qff \not\in\qff l_{\fff i}$\oss
for all\qss 
$i\off =\off
1\fff,\pff 2\fff,\pff \ldots\fff,\pff p$\dnsp,\oss 
i.e.\qss on the fact that\oss
$i\qff \longmapsto\qff a_{\fff i\dff +\dff 1}$\oss
is a system of distinct representatives for the family\oss
$i\off \longmapsto\off E\qff \smallsetminus\qff l_{\fff i}$\oss
of subsets of $E$\nnsp,\oss
where\oss
$i\off =\off
1\fff,\pff 2\fff,\pff \ldots\fff,\pff p$\dnsp.\oss

Once this is realized,\oss
it is only natural to look for a system of distinct representatives of
the full family\oss 
$\dis
l\qff \longmapsto\qff E\qff \smallsetminus\qff l$\oss
of the complements of lines,\oss
i.e.\qss for an injective map\oss
$l\qff \longmapsto\qff a\fff({\fff}l\fff)$\oss from\dss $\mathcal{L}$\dss to\dss $E$\dss 
such that\oss
$a\fff({\fff}l\fff)\off \in\off E\qff \smallsetminus\qff l$\oss
for all\qss $l\qff \in\qff \mathcal{L}$\dnsp.\oss

By the well known\dss Ph.\dss Hall's marriage theorem,\oss
such a system of distinct representatives exists if and only if
for every subset\qss $\mathcal{K}\qff \subset\qff \mathcal{L}$\qss
the union\vspace*{3pt}
\begin{equation}
\label{union}
\quad
\bigcup_{\dff l\dff \in\dff \mathcal{K}}\qff E\qff \smallsetminus\qff l
\off\off =\off\off
E\off \smallsetminus\off \bigcap_{\dff l\dff \in\dff \mathcal{K}}\qff l
\end{equation}

\vspace*{-9pt}
contains\qss $\geqslant\qff |\fff \mathcal{K}\fff |$\qss elements,\oss
where\dss $|\dff X\dff|$\dss denotes the number of elements of a set $X$\nnsp.\oss
But the intersection of\qss $\geqslant\qff 2$\qss lines consists of\qss
$\leqslant\qff 1$\qss points,\oss and,\oss
almost obviously,\oss this condition holds.

\myuppar{The message.} 
All this emerged in my mind in one instant as an irreducible revelation.\oss
My first thought after this revelation
was that it cannot be true,\oss because if it is true,\oss
then everybody writing about this topic would 
use systems of distinct representatives.\oss
Perhaps,\pss the right question is not how I came up with this idea,\pss
but why experts missed it\halfff.
The rest of this section is devoted to the proof\qss \cite{i-db-e}\qss
based on this revelation.

\myuppar{Proof\dss of\qss $m\qff \geqslant\qff n$\nnsp.}
We may assume that $m\qff \leqslant\qff n$\nnsp.\oss
Let $\mathcal{K}$ be a subset of $\mathcal{L}$\dnsp.\oss
If\qss $|\dff \mathcal{K}\dff|\qff =\qff 1$\nnsp,\qss
then\qss (\ref{union})\qss is the complement of a line 
and hence contains\qss $\geqslant\qff 1$\qss elements.\qss 
If\qss $2\qff \leqslant\qff |\dff \mathcal{K}\dff |\qff \leqslant\qff m\qff -\qff 1$\nnsp,\qss
then\qss (\ref{union})\qss 
is a complement in $E$ of\qss $\leqslant\qff 1$\qss point
and hence contains\qss
\[
\quad
\geqslant\qff n\qff -\qff 1
\off \geqslant\off 
m\qff -\qff 1
\off \geqslant\off
|\dff \mathcal{K}\dff |
\] 
elements.\qss
If\qss $|\dff \mathcal{K}\dff|\qff =\qff m$\nnsp,\qss
then\qss (\ref{union})\qss contains\qss
$n\qff \geqslant\qff m\qff =\qff |\dff \mathcal{K}\dff|$\qss elements.\qss
Therefore there exists a system of distinct representatives\oss
for the family\oss 
$l\qff \longmapsto\qff E\qff \smallsetminus\qff l$\dnsp,\oss
i.e. there exists an injective map\oss
$l\qff \longmapsto\qff a\fff({\fff}l\fff)$\oss such that\qss
$a\fff({\fff}l\fff)\qff \not\in\qff l$\qss for every  $l$\dnsp.\oss
By the \dbe inequalities
\begin{equation}
\label{basic-reps-ineq}
\quad
s_{\fff l}\off \leqslant\off k_{\fff a\fff({\fff}l\fff)}
\quad\mbox{ for\qss every }\quad
l\qff \in\qff \mathcal{L}.
\end{equation}
By summing all these inequalities and using the injectivity of\oss
$l\qff \longmapsto\qff a\fff({\fff}l\fff)$\oss 
we see that
\begin{equation}
\label{three-sums}
\quad
\sum_{l\dff \in\dff \mathcal{L}}\qff s_{\fff l}
\off\off \leqslant\off\off
\sum_{l\dff \in\dff \mathcal{L}}\qff k_{\fff a\fff({\fff}l\fff)}
\off\off \leqslant\off\off
\sum_{\qff z\dff \in\dff E}\qff k_{\fff z}\qff.
\end{equation}
Moreover\halfff,\qss
the second inequality is strict unless\qss $m\qff =\qff n$\qss
(otherwise the last sum has more positive summands
than the previous one).\oss
But\qss (\ref{sums})\qss implies that both inequalities in\qss
(\ref{three-sums})\qss should be actually equalities.\oss
It follows that\qss $m\qff =\qff n$\nnsp.\oss
Moreover\halfff,\oss in view of the inequalities\qss
(\ref{basic-reps-ineq}),\oss
it follows that\oss 
$s_{\fff l}\off =\off k_{\fff a\fff({\fff}l\fff)}$\oss
for every\qss $l\qff \in\qff \mathcal{L}$\qss
(under the assumption\qss $m\qff \leqslant\qff n$\nnsp).

\myuppar{The case\qss $m\qff =\qff n$\nnsp.}
Suppose that a point $z$
is contained in\qss $\geqslant\qff m\qff -\qff 1$\qss lines.\oss
Each of these lines contains at least one point in addition to $z$\nnsp.\oss
Since\qss $m\qff =\qff n$\nnsp,\qss
there are no other points and $z$ is contained in exactly\qss $m\qff -\qff 1$\qss lines.\oss
Since there are exactly $m$ lines,\oss
only one line does not contain $z$\nnsp.\oss
This line should contain all points\qss $\neq\qff z$\nnsp.\qss
It follows that\qss $(\dff E\fff,\pff \mathcal{L}\dff)$\qss is a near-pencil.\oss

Suppose now that no point 
is contained in\qss $\geqslant\qff m\qff -\qff 1$\qss lines.\oss
Let $\mathcal{K}$ be a proper subset of $\mathcal{L}$\dnsp.\oss
If\qss $|\dff \mathcal{K}\dff|\qff =\qff 1$\nnsp,\oss 
then\qss (\ref{union})\qss is equal to\qss $E\dff \smallsetminus\dff l$\qss
for some line $l$\nnsp.\oss
If\qss $E\dff \smallsetminus\dff l$\qss consists of only one point $z$\nnsp,\oss
then by the \dbe inequalities $z$ is contained in\qss
$\geqslant\qff s_{\fff l}\qff =\qff m\qff -\qff 1$\qss lines,\oss
contrary to the assumption.\oss
Therefore,\oss (\ref{union})\qss contains\qss
$\geqslant\qff 2\qff =\qff |\dff \mathcal{K}\dff|\qff +\qff 1$\qss points.\oss
If\qss  $|\dff \mathcal{K}\dff|\qff \leqslant\qff m\qff -\qff 2$\nnsp,\oss
then\qss (\ref{union})\qss contains\qss
$\geqslant\qff n\qff -\qff 1\qff =\qff m\qff -\qff 1\qff \geqslant\qff
|\dff \mathcal{K}\dff|\qff +\qff 1$\qss points.\oss
Finally,\oss if\qss  $|\dff \mathcal{K}\dff|\qff =\qff m\qff -\qff 1$\nnsp,\oss
then\qss (\ref{union})\qss contains all\qss 
$n\qff =\qff m\qff =\qff |\dff \mathcal{K}\dff|\qff +\qff 1$\qss 
points because no point 
is contained in\qss $\geqslant\qff m\qff -\qff 1$\qss lines.\oss

We see that\qss (\ref{union})\qss contains\qss
$\geqslant\qff |\dff \mathcal{K}\dff|\qff +\qff 1$\qss elements
for every proper subset\qss $\mathcal{K}\qff \subset\qff \mathcal{L}$\dnsp.\oss 
This allows to get from the marriage theorem more than
just the existence of a system of distinct representatives.\oss
Let\qss $\lambda\qff \in\qff \mathcal{L}$\qss and\qss 
$z\qff \in\qff E\dff \smallsetminus\qff \lambda$\nnsp.\oss
Then there exists a system of distinct representatives\oss
$l\qff \longmapsto\qff a\fff({\fff}l\fff)$\oss such that\qss
$a\fff(\fff\lambda\fff)\qff =\qff z$\nnsp.\oss
This immediately follows from an application
of the marriage theorem
to the family of sets\oss
$(\dff E\qff \smallsetminus\qff \{\dff z \qff\}\dff)\qff \smallsetminus\qff l$\oss
with\oss
$l\qff \in\qff \mathcal{L}\qff \smallsetminus\qff \{\dff \lambda \qff\}$\nnsp.\oss

Since\qss
$m\qff \leqslant\qff n$\nnsp,\oss
the existence of a system of distinct representatives\oss
$l\qff \longmapsto\qff a\fff({\fff}l\fff)$\oss such that\qss
$a\fff(\fff\lambda\fff)\qff =\qff z$\qss
implies that\oss
$s_{\fff \lambda}\off =\off k_{\fff a\fff({\fff}\lambda\fff)}\off =\off k_{\fff z}$\nnsp.\oss
Therefore,\qss $z\qff \not\in\qff l$\qss implies that\qss
$s_{\fff l}\qff =\qff k_{\fff z}$\pss
and hence every line containing $z$ intersects $l$\dnsp.\oss
It follows that every two lines intersect\halfff.\oss

If $E$ cannot be obtained as the union of two lines,\oss
then for every two lines\qss $l\fff,\pff l'$\qss
there exists a point $z$ such that\qss $z\qff \not\in\qff l\fff,\pff l'$\qss
and hence\qss $s_{\fff l}\qff =\qff k_{\fff z}\qff =\qff s_{\fff l'}$\nnsp.\oss
In this case all the numbers\qss 
$s_{\fff l}\fff,\pff k_{\fff z}$\qss are equal
and hence\qss $(\dff E\fff,\pff \mathcal{L}\dff)$\qss is a projective plane
by Lemma\qss 1\qss at the end of Section\qss \ref{solution}.\oss
If there exist two lines\qss $l\fff,\pff l'$\qss
such that\pss $E\qff =\qff l\dff \cup\dff l'$\nnsp,\qff\oss
then\qss $k_{\fff y}\qff =\qff 2$\nnsp,\oss
where $y$ is the point of intersection of $l$ and $l'$\dnsp,\oss
and the proof\dss is completed by applying the following lemma.

\myuppar{Lemma\qss 2.}
\emph{If\qss $m\qff =\qff n$\qss and\qss $k_{\fff y}\qff =\qff 2$\qss for some point $y$\dnsp,\oss
then\qss $(\dff E\fff,\pff \mathcal{L}\dff)$\qss is a near-pencil.}

\myuppar{Proof\halfff.}
Let\qss $l\fff,\pff l'$\qss be the lines containing $y$\dnsp.\oss
Then\qss $E\qff =\qff l\dff \cup\dff l'$\qss and there are\qss
$n\qff =\qff s_{\fff l}\qff +\qff s_{\fff l'}\qff -\qff 1$\qss points.\oss
In addition to the lines\pss $l\fff,\pff l'$\pss
there\dss are\qss
$(s_l\qff -\qff 1)(s_{l'}\qff -\qff 1)$\qss
lines consisting of a point in\qss 
$l\qff \smallsetminus\dss \{\dff y \trf\}$\qss
and a point in\qss 
$l'\qff \smallsetminus\dss \{\dff y \trf\}$\nnsp.\qff\oss
If\oss $s_{\fff l}\off \geqslant\off s_{\fff l'}\off \geqslant\off 3$\dnsp,\oss
then the number $m$\dss of\trs lines\dss is\qss\vspace*{1.5pt}
\[
\quad
\geqslant\off
2\qff +\qff (s_l\qff -\qff 1)(s_{l'}\qff -\qff 1)
\off \geqslant\off 
2\qff +\qff 2\dff(s_l\qff -\qff 1)
\off =\off 
2\fff s_l
\off \geqslant\off 
s_l\qff +\qff s_{l'}
\off =\off 
n\qff +\qff 1\dff,
\]

\vspace*{-10.5pt}
contrary to the assumption\qss $m\qff =\qff n$\nnsp.\qss
Therefore one of the lines\qss $l\fff,\pff l'$\qss
consists of $2$ points
and hence\qss $(\dff E\fff,\pff \mathcal{L}\dff)$\qss 
is a near-pencil.\oss  \eproof

\mysection{Linear\qss algebra\qss and\qss the\qss inequality\qss
$m\qff \geqslant\qff n$}{linear-algebra}

\vspace*{6pt}
\myuppar{A proof of the inequality\qss
$m\qff \geqslant\qff n$\qss based on the linear independence.}
This proof was communicated to me by\qss F.\qss Petrov\qss \cite{p}.\pss
I\dss believe that this is essentially the proof found by A.\qss Suslin.

Let\qss $\rl$\qss be the vector space of maps\qss $\mathcal{L}\dff \toto\dff \mathbf{R}$\qss
with the scalar product
\[
\quad
(\dff v\fff,\pff w \dff)
\off\off =\off\off
\sum_{l\qff \in\qff \mathcal{L}}\qff v(\fff l\fff)\dff w(\fff l\fff).
\]
Every\qss $z\qff \in\qff E$ defines a map\qss 
$v_{\fff z}\dff \colon\dff \mathcal{L}\dff \toto\dff \mathbf{R}$\qss
by the rule\qss $v_{\fff z}\fff(\fff l\fff)\qff =\qff 1$\qss
if\qss $z\qff \in\qff l$\qss and\qss
$v_{\fff z}\fff(\fff l\fff)\qff =\qff 0$\qss otherwise.\qss
There are $n$ maps $v_{\fff z}$\nnsp.\oss
Since the dimension of $\rl$ is equal to $m$\nnsp,\oss
it is sufficient to prove that the maps $v_{\fff z}$ 
are independent as vectors of\dss $\rl$\dnsp.\oss

The scalar product\qss $(\dff v_{\fff z}\fff,\pff v_{\fff z} \dff)$\qss
is equal to the number of lines containing the point $z$\nnsp,\oss
and hence\qss
$(\dff v_{\fff z}\fff,\pff v_{\fff z} \dff)
\qff \geqslant\qff 2$\qss
for all\qss $z\qff \in\qff E$\nnsp.\oss
If\qss $z\qff \neq\qff y$\dnsp,\oss then\qss 
$(\dff v_{\fff z}\fff,\pff v_{\fff y} \dff)$\qss
is equal to the number of lines containing both $z$ and $y$\nnsp,\pss
and hence\qss
$(\dff v_{\fff z}\fff,\pff v_{\fff y} \dff)
\qff =\qff 1$\nnsp.\oss
If the vectors $v_{\fff z}$ 
are linearly dependent\halfff,\oss
then
\begin{equation*}
\quad
\sum_{z\qff \in\qff E}\qff c_{\fff z}\dff v_{\fff z}
\off\off =\off\off
0
\end{equation*}
for some real numbers\qss $c_{\fff z}$\nnsp,\qss $z\qff \in\qff E$\nnsp,\oss
such that not all $c_{\fff z}$ are equal to $0$\nnsp.\oss
For every\qss $y\qff \in\qff E$\qss
taking the scalar product of this equality with the vector\dss $v_{\fff y}$\dss
results in the equality
\[
\quad
\sum_{z\qff \in\qff E}\qff c_{\fff z}\dff (\fff v_{\fff z}\fff,\pff v_{\fff y} \fff)
\off\off =\off\off
0.
\]
Since\qss $(\fff v_{\fff z}\fff,\pff v_{\fff y} \fff)\qff =\qff 1$\qss
for all\qss $z\qff \neq\qff y$\nnsp,\oss
this equality implies that
\[
\quad
c_{\fff y}\dff ((\fff v_{\fff y}\fff,\pff v_{\fff y} \fff)\qff -\qff 1)
\off +\off
\sum_{z\qff \in\qff E}\qff c_{\fff z} 
\off\off =\off\off
0.
\]
Since\qss $(\fff v_{\fff y}\fff,\pff v_{\fff y} \fff)\qff \geqslant\qff 2$\nnsp,\oss
it follows that the coefficient $c_{\fff y}$ and the sum\oss
\[
\quad
\sum_{z\qff \in\qff E}\qff c_{\fff z}
\] 
have opposite signs.\oss
But since not all $c_{\fff y}$ are equal to $0$\nnsp,\oss 
this cannot be true for all\qss $y\qff \in\qff E$\nnsp.\oss
The contradiction shows that vectors $v_{\fff z}$ are linearly independent
and hence\qss $m\qff \geqslant\qff n$\nnsp.\oss  \eproof

\myuppar{Standard linear algebra proofs of the inequality\qss
$m\qff \geqslant\qff n$\nnsp.}
In order to present standard proofs it is convenient to return to
the notations of the\dss N.\qss Bourbaki exercise.\oss
Let $M$ be the\qss \emph{incidence\dss matrix}\qss
of the points $a_j$ and sets $A_i$\nnsp.\oss
Namely,\qss $M$\dss is\dss an $n \times m$ matrix with entries\qss
$m_{j i}\qff =\qff 1$\qss if\qss $a_j\qff \in\qff A_i$\qss
and\qss $m_{j i}\qff =\qff 0$\qss otherwise.\oss
Let us consider the product $M\fff M^{\fff T}$\nnsp,\oss
where $M^{\fff T}$ is the matrix transposed to $M$\nnsp.\oss
It is an $n \times n$ matrix with all non-diagonal entries equal to $1$
and with diagonal entries\oss
$k_{\fff 1}\fff,\pff  k_{\fff 2}\fff,\pff \ldots\fff,\pff k_{\fff n}$\nnsp.\oss
The most classical linear algebra proofs\halfff,\oss
going back to the paper\qss \cite{bo}\qss by\qss R.C.\qss Bose,\oss 
proceed with the computation of the determinant of
$M\fff M^{\fff T}$\dnsp.\oss
It is rarely presented in details;\oss
apparently,\pss it is expected that the readers enjoy computations of determinants.\oss
Curious readers may find a computation of\qss $\det\qff M\fff M^{\fff T}$\qss
at the end of this section;\oss
in particular\halfff,\oss
the computation shows that this determinant is non-zero.\oss
The non-vanishing of\qss $\det\qff M\fff M^{\fff T}$\qss
means that the rank of the matrix $M\fff M^{\fff T}$ is equal to $n$\nnsp,\oss
and this implies that the rank of $M$ is\qss $\geqslant\qff n$\nnsp.\oss
Since $M$ is\dss an $n \times m$ matrix,\oss
this may happen only if\qss $m\qff \geqslant\qff n$\nnsp.\oss

More modern expositions avoid computation of the determinant\qss
$\det\qff M\fff M^{\fff T}$\qss
by observing that $M\fff M^{\fff T}$ is equal to the sum of the diagonal matrix
with the diagonal entries\qss\vspace*{2pt}
\[
\quad
k_{\fff 1}\qff -\qff 1\fff,\off  k_{\fff 2}\qff -\qff 1\fff,\off \ldots\fff,\off k_{\fff n}\qff -\qff 1
\]

\vspace*{-10pt}
and the $n \times n$ matrix $J$ with all entries equal to $1$\nnsp.\oss
Since\qss $k_{\fff j}\qff \geqslant\qff 2$\qss and hence\qss
$k_{\fff j}\qff -\qff 1\qff \geqslant\qff 1$\qss for all $j$\dnsp,\oss
the first matrix is positive definite.\oss
The matrix $J$ is positive semi-definite,\qss
although is not definite.\oss
In fact\halfff,\oss
the associated quadratic form $\mathbf{x}\qff J\trf \mathbf{x}^{\fff T}$\dnsp,\oss
where\qss
$\mathbf{x}
\off =\off
(\fff x_{\fff 1}\fff,\pff x_{\fff 2}\fff,\pff \ldots,\fff x_{\fff n} \fff)$\qss
is a row vector\halfff,\oss 
is equal to\qss
$(\fff x_{\fff 1}\qff +\qff x_{\fff 2}\qff +\qff \ldots\qff +\qff x_{\fff n} \fff)^{\fff 2}$\dnsp.\oss
It follows that the sum $M\fff M^{\fff T}$ of these matrices
is positive definite and hence has the rank $n$\nnsp.\oss
As above,\oss this implies that\qss $m\qff \geqslant\qff n$\nnsp.\oss

\myuppar{Comparing the proofs.}
The standard proofs do not fit the\qss 
\emph{``Kvant''}\qss description of the proof by A.\qss Suslin:\qss
they use more advanced tools than the theorem about the linear dependence of
more than $n$ vectors in an $n$\dnsp-dimensional vector space.\oss
One can find a proof based only on this theorem
in the unpublished book draft\qss \cite{bf}\qss by\dss
L.\qss Babai and\dss P.\qss Frankl.\oss
But even in this remarkable book it is hidden in the exercises.\oss
See Exercise\qss 4.1.5\qss and its solution on p.\qss 184.\oss
The preference for using the matrix $M\fff M^{\fff T}$
seems to be a part of a dominating culture.\oss
On the other hand,\oss all proofs based on the linear algebra more or less explicitly reduce
the inequality\qss $m\qff \geqslant\qff n$\qss to the following lemma and then prove it\halfff.\oss

\myuppar{Lemma.}
\emph{Let\qss $V$ be an $m$\dnsp-dimensional vector space over\qss $\rr$ equipped with
a scalar product\qss $(\fff \bullet\fff,\pff \bullet \fff)$\dnsp.\oss
Let\qss $P$ be a set of $n$ vectors in\qss $V$\dnsp.\oss
Suppose that there exists\qss 
$\lambda\qff \in\qff \rr$\nnsp,\oss
$\lambda\qff >\qff 0$\dnsp,\oss
such that\qss}\vspace*{2pt} 
\[
\quad
(\fff u\fff,\pff u \fff)\off >\off \lambda
\hspace*{1.5em}\mbox{ \emph{and} }\hspace*{1.5em}
(\fff v\fff,\pff w \fff)\off =\off \lambda
\]

\vspace*{-10pt}
\emph{for every\qss $u\qff \in\qff P$\qss
and every two distinct vectors\qss
$v\fff,\pff w\qff \in\qff P$\dnsp.\oss
Then\qss $m\qff \geqslant\qff n$\nnsp.\oss}  \eproof

\myuppar{A generalization.}
The linear algebra proofs apply with only trivial changes to a more general situation.\oss
Namely,\oss it is sufficient to assume that 
there exist a natural number\qss $\lambda\qff \geqslant\qff 1$\qss
such that every two distinct points are contained in exactly $\lambda$ lines
and every point is contained in\qss $>\qff \lambda$\qss lines.\oss
Then the conclusion\qss $m\qff \geqslant\qff n$\qss still holds.\oss
This is due to H.J.\dss Ryser\qss \cite{r}.\oss
Apparently,\oss no combinatorial proof of\dss
Ryser's theorem is known.\oss
Ryser\qss \cite{r}\qss also used linear algebra to provide a description of
the case\qss $m\qff =\qff n$\qss similar to
de Bruijn--Erd\"{o}s description in the case\qss $\lambda\qff =\qff 1$\nnsp.\oss

\myuppar{The determinant of\qss $M\fff M^{\fff T}$\nnsp.}
For the benefit of the readers who do not like 
to compute the determinants themselves\halfff,\oss 
here is a computation of\qss 
$\det\qff M\fff M^{\fff T}$\qss 
following the textbook\qss \cite{hp}.\oss

Let\qss $m_{j}\qff =\qff k_{\fff j}\qff -\qff 1$\qss
for all $j$\dnsp.\oss
Then\vspace*{6pt}
\[
\quad
M\fff M^{\fff T}
\off =\off
\begin{bmatrix}\off
m_{\fff 1}\dff +\qff 1 & 1 & 1 & \hdotsfor{4} & 1 & 1\\
1 & m_{\fff 2}\dff +\qff 1 & 1 & \hdotsfor{4} & 1 & 1\\
1 & 1 & m_{\fff 3}\dff +\qff 1 & \hdotsfor{4} & 1 & 1\\
& & & & & & & &  \\
 \hdotsfor{9}    \\
& & & & & & & &  \\
1 & 1 & 1 & \hdotsfor{4} & 1 & m_{\fff n}\dff +\qff 1\off
\end{bmatrix}.
\]

\vspace*{-6pt}
Let us subtract the first row from every other one and get the matrix\vspace*{6pt}
\[
\quad
\phantom{M\fff M^{\fff T}
\off =\off}
\begin{bmatrix}\off
m_{\fff 1}\dff +\qff 1 & 1 & 1 & \hdotsfor{4} & 1 & 1\\
-\qff m_{\fff 1}  & m_{\fff 2} & 0 & \hdotsfor{4} & 0 & 0\\
-\qff m_{\fff 1} & 0 & m_{\fff 3} & \hdotsfor{4} & 0 & 0\\
& & & & & & & &  \\
 \hdotsfor{9}    \\
& & & & & & & &  \\
-\qff m_{\fff 1} & 0 & 0 & \hdotsfor{4} & 0 & m_{\fff n}\off
\end{bmatrix}.
\]

\vspace{-6pt}
For\oss
$j\off =\off 2\fff,\pff 3\fff,\pff \ldots\fff,\pff n$\nnsp,\oss
let us multiply the $j$\dnsp-th column by\qss $m_{\fff 1}/m_{\fff j}$\qss
(recall that\qss $k_{\fff j}\qff \geqslant\qff 2$\qss 
and hence\qss $m_j\qff \geqslant\qff 1\qff >\qff 0$\nnsp)\qss
and add the result to the first column.\oss
We get the matrix\vspace*{6pt}
\[
\quad
\phantom{M\fff M^{\fff T}
\off =\off}
\begin{bmatrix}\off
D & 1 & 1 & \hdotsfor{4} & 1 & 1\\
0  & m_{\fff 2}  & 0 & \hdotsfor{4} & 0 & 0\\
0 & 0 & m_{\fff 3} & \hdotsfor{4} & 0 & 0\\
& & & & & & & &  \\
 \hdotsfor{9}    \\
& & & & & & & &  \\
0 & 0 & 0 & \hdotsfor{4} & 0 & m_{\fff n}\off
\end{bmatrix},
\]

\vspace*{-4pt}
where\off\oss
$\dis
D
\off\off =\off\off
m_{\fff 1}
\qff +\qff
1
\off +\off
\sum_{j\qff =\qff 2}^n\qff \frac{\fff m_{\fff 1}\fff}{m_j}
\off\off =\off\off
m_{\fff 1}
\off +\off
\sum_{j\qff =\qff 1}^n\qff \frac{\fff m_{\fff 1}\fff}{m_j}$\qff.\vspace*{2pt}

It follows that\off\oss 
$\dis
\det\qff M\dff M^{\fff T}
\off\off =\off\off
D
\qff\cdot 
\prod_{j\qff =\qff 2}^n m_j
\off\off =\off\off
\prod_{j\qff =\qff 1}^n m_j
\qff \cdot\qff
\left(\qff
1
\qff +\qff
\sum_{j\qff =\qff 1}^n\qff \frac{\fff 1\fff}{m_j}
\qff\right)
\off\off \neq\off\off 0$.\oss

\mysection{Hanani's\qss theorem}{h-proof}

\vspace*{6pt}
\myuppar{Two papers of\dss H.\dss Hanani.}
According to the Th.\dss Motzkin\qss \cite{m},\oss
the first proof of the inequality\dss $m\qff \geqslant\qff n$\dss
and,\pss it seems,\pss
of the full de Bruijn--Erd\"{o}s theorem,\oss
was given in\qss 1938\qss by H.\dss Hanani.\oss
He published an outline of his proof\qss \cite{h1}\qss 
only in\qss 1951.\oss 
Later on H.\dss Hanani published a detailed 
exposition\qss \cite{h2}\qss of a simplified proof\halfff.\oss 
In fact\halfff,\oss in\qss \cite{h2}\qss he proved\qss
(at no extra cost)\qss
a stronger version of the
de Bruijn--Erd\"{o}s theorem.\oss
He also used his methods to prove a\dss $3$\dnsp-dimensional analogue
dealing with points,\qss lines,\qss and planes.\oss

\myuppar{Hanani's\dss Theorem.}
\emph{Under the previous assumptions,\oss
let\qss $L\qff \in\qff \mathcal{L}$\qss be a line containing the maximal
number of points among all lines,\oss
let\qss $\mathcal{P}$\qss be the set of all lines intersecting\qss $L$\qss
(in particular\halfff,\oss $L\qff \in\qff \mathcal{P}$\nnsp),\oss
and let\dss $p$\dss be the number of elements of\qss $\mathcal{P}$\dnsp.\oss
Then\qss $p\qff \geqslant\qff n$\nnsp,\oss
and\dss if\qss $p\qff =\qff n$\nnsp,\oss then\qss
$\mathcal{P}\qff =\qff \mathcal{L}$\qss and\pss $(\dff E\fff,\pff \mathcal{L}\dff)$\pss 
is either a near-pencil,\oss or a projective plane.\oss}

\vspace*{6pt}
Suppose
that\dss $n\qff \geqslant\qff p$\dss
and\qss $(\dff E\fff,\pff \mathcal{L}\dff)$\qss is not a near-pencil.\oss
As usual,\oss we assume that every line contains\qss $\geqslant\qff 2$\qss points.\oss
Let\qss $a\qff =\qff s_{\dff L}$\nnsp.\oss
Let $K$ be the line with the maximal number of points among the lines different from $L$\nnsp,\oss
and\dss let\qss $b\qff =\qff s_{\dff K}$\nnsp.\oss
Then\qss $a\qff \geqslant\qff b$\nnsp.\oss
The strategy is to estimate $n$\nnsp,\pss
or\halfff,\pss what is the same,\qss $n\qff -\qff 1$\qss
in terms of $a$ and $b$ both from the below and from the above.\oss

\myuppar{Hanani's\dss Lemma.}
\emph{If\qss $x\qff \in\qff L$\nnsp,\oss then}\oss\vspace*{2pt}
\begin{equation}
\label{lemma}
k_{\fff x}\qff -\qff 1\off \geqslant\off \frac{n\qff -\qff a}{b\qff -\qff 1}\qff.
\end{equation}

\vspace*{-10pt}
\proof\qss
Let us consider pairs\qss $(\fff l\fff,\pff y\fff)$\qss
such that $l$ is a line containing $x$ and 
$y$ is a point in\qss $l\dff \smallsetminus\dff L$\nnsp.\oss
Such a pair is uniquely determined by the point $y$ 
and hence there are $n\qff -\qff a$ such pairs.\oss
A line $l$ occurs in such a pair if and only if\qss $x\qff \in\qff l$\qss
and\qss $l\qff \neq\qff L$\nnsp.\oss
It follows that there are\qss $k_{\fff x}\qff -\qff 1$\qss choices of $l$\nnsp.\oss
Given a line $l$\dnsp,\oss 
there are\qss $\leqslant\qff b\qff -\qff 1$\qss
choices for the point $y$\dnsp.\oss
Therefore the number\qss $n\qff -\qff a$\qss of such pairs is\qss 
$\leqslant\qff (k_{\fff x}\qff -\qff 1)(b\qff -\qff 1)$\dnsp.\oss
The lemma follows.\oss  \eproof

\myuppar{An upper estimate of\qss $n\qff -\qff 1$\nnsp.}
By summing the inequalities\qss (\ref{lemma})\qss over all\qss $x\qff \in\qff L$\qss
and adding $1$ in order to account for the line $L$ itself\halfff,\oss
we can estimate $p$ from below and conclude that\vspace*{3pt}
\begin{equation}
\label{p-lower-1}
\quad
n
\off \geqslant\off
p
\off \geqslant\off
1\qff +\qff a\qff \frac{n\qff -\qff a}{b\qff -\qff 1}
\off =\off
1\qff +\qff a\qff \frac{(n\qff -\qff 1)\qff -\qff (a\qff -\qff 1)}{b\qff -\qff 1}
\end{equation}

\vspace*{-9pt}
or\halfff,\oss what is the same,\oss\vspace*{3pt}
\begin{equation}
\label{n-upper}
\quad
a\fff(a\qff -\qff 1)
\off \geqslant\off
(n\qff -\qff 1)(a\qff -\qff b\qff +\qff 1)\dff.
\end{equation}

\vspace*{-9pt}
The inequality\qss (\ref{n-upper})\qss 
provides an estimate of\qss $n\qff -\qff 1$\qss from the above.\oss

\myuppar{A lower estimate of\qss $n\qff -\qff 1$\nnsp.}
There is another way to estimate $p$ from below.\oss
By a miracle,\oss this other estimate of $p$ from the same side 
leads to an estimate of\qss $n\qff -\qff 1$\qss from the
other side.\oss
Let $z$ be a point in\qss $L\dff \cap\dff K$\qss if\qss 
$L\dff \cap\dff K\qff \neq\qff \emptyset$\dnsp,\oss
and an arbitrary point of $L$ otherwise.\oss
For every\qss $x\qff \in\qff L\dff \smallsetminus\dff \{\dff z\dff\}$\nnsp,\qss
$y\qff \in\qff K\dff \smallsetminus\dff \{\dff z\dff\}$\qss
there is a unique line $l$ containing\qss $\{\dff x\fff,\pff y \dff\}$\nnsp.\oss
All these lines are distinct\halfff,\oss not equal to $L$\nnsp,\oss 
and do not contain $z$\nnsp.\oss 
Clearly,\oss there are\qss
$(a\qff -\qff 1)(b\qff -\qff 1)$\qss
of such lines.\oss
A lower estimate of number $k_{\fff z}$
of lines containing $z$ is provided by\qss (\ref{lemma}).\oss
It follows that\vspace*{4pt}
\begin{equation*}
\quad
n
\off \geqslant\off
p
\off \geqslant\off
1
\qff +\qff 
\frac{n\qff -\qff a}{b\qff -\qff 1}
\qff +\qff
(a\qff -\qff 1)(b\qff -\qff 1)
\end{equation*}

\vspace*{-33pt}
\[
\quad
\hspace*{12em}
\phantom{n
\off \geqslant\off
p
\off }
=\off
1
\qff +\qff 
\frac{(n\qff -\qff 1)\qff -\qff (a\qff -\qff 1)}{b\qff -\qff 1}
\qff +\qff
(a\qff -\qff 1)(b\qff -\qff 1)
\]

\vspace*{-8pt}
and hence\qff\oss
$\dis
(n\qff -\qff 1)(b\qff -\qff 1)
\off \geqslant\off
(n\qff -\qff 1)
\qff -\qff
(a\qff -\qff 1)
\qff +\qff
(a\qff -\qff 1)(b\qff -\qff 1)^{\fff 2}$\qff\oss
and\vspace*{3pt}
\begin{equation*}
\quad
(n\qff -\qff 1)(b\qff -\qff 2)
\off \geqslant\off
(a\qff -\qff 1)(b^{\fff 2}\qff -\qff 2\fff b)
\off =\off
(a\qff -\qff 1)\fff b\fff (b\qff -\qff 2).
\end{equation*}

\vspace*{-9pt}
Since\qss $b\qff \geqslant\qff 2$\nnsp,\oss
it follows that either\qss $b\qff =\qff 2$\nnsp,\oss
or\qss
$\dis
n\qff -\qff 1
\off \geqslant\off 
(a\qff -\qff 1)\fff b$\nnsp.\oss
If\qss $b\qff =\qff 2$\nnsp,\oss 
then all lines except $L$ consist of $2$ points and the inequality\qss 
(\ref{n-upper})\qss implies that\qss
$a\qff \geqslant\qff n\qff -\qff 1$\nnsp.\oss
But\qss $L\qff \neq\qff E$\qss and hence\qss
$a\qff \leqslant\qff n\qff -\qff 1$\nnsp.\oss
It follows that\qss $a\qff =\qff n\qff -\qff 1$\qss
and hence $L$ contains all points of $E$ except one
and\qss $(\dff E\fff,\pff \mathcal{L}\dff)$\qss is a near-pencil,\oss
contrary to the assumption.\oss
Therefore\qss\vspace*{3pt}
\begin{equation}
\label{n-lower}
\quad
n\qff -\qff 1
\off \geqslant\off 
(a\qff -\qff 1)\fff b\dff.
\end{equation}

\vspace*{-9pt}
The inequality\qss (\ref{n-lower})\qss 
provides an estimate of\qss $n\qff -\qff 1$\qss from the below.\oss

\myuppar{Combining the two estimates.}
After multiplying the inequality\qss (\ref{n-lower})\qss by\qss
$(a\qff -\qff b\qff +\qff 1)$\qss
and combining the result with the inequality\qss (\ref{n-upper}),\oss 
we see that\vspace*{3pt}
\[
\quad
a\fff (a\qff -\qff 1)
\off \geqslant\off 
(a\qff -\qff 1)\fff b\fff (a\qff -\qff b\qff +\qff 1)
\]

\vspace*{-9pt}
and hence\oss
$\dis
a
\off \geqslant\off
b\fff(a\qff -\qff b\qff +\qff 1)
\off =\off
b\fff (a\qff -\qff b)
\qff +\qff
b$\nnsp,\oss
or\halfff,\oss what is the same
\[
\quad
0
\off \geqslant\off
(b\qff -\qff 1) (a\qff -\qff b)\dff.
\]
Since\qss $b\qff >\qff 1$\nnsp,\oss this implies that\qss $b\qff \geqslant\qff a$\nnsp.\oss
On the other hand,\oss
$b\qff \leqslant\qff a$\qss by the definition of\qss $a\fff,\pff b$\nnsp.\oss
It follows that\qss $a\qff =\qff b$\nnsp.\oss
By combining\qss $a\qff =\qff b$\qss with the inequalities\qss
(\ref{n-upper})\qss
and\qss (\ref{n-lower})\qss
we conclude,\oss respectively,\oss
that\qss 
$a\fff(a\qff -\qff 1)
\off \geqslant\off
n\qff -\qff 1$\qss
and\qss
$n\qff -\qff 1
\off \geqslant\off
a\fff(a\qff -\qff 1)$\dnsp.\oss
It follows that\qss
$n\qff -\qff 1
\off =\off
a\fff(a\qff -\qff 1)$\qss
and\dss hence\qss
$n
\off =\off
a\fff(a\qff -\qff 1)\qff +\qff 1$\nnsp.\oss
By combining this with\qss (\ref{p-lower-1})\qss we see that\vspace*{3pt}
\[
\quad
a\fff(a\qff -\qff 1)\qff +\qff 1
\off =\off
n
\off \geqslant\off
p
\off \geqslant\off
1\qff +\qff a\qff \frac{a\fff(a\qff -\qff 1)\qff +\qff 1\qff -\qff a}{a\qff -\qff 1}
\off =\off
1\qff +\qff a\fff(a\qff -\qff 1)\dff.
\]

\vspace{-9pt}
It follows that\qss $n\qff =\qff p$\dnsp,\oss
and therefore\qss $p\qff \geqslant\qff n$\qss
if the inequality\qss $p\qff \leqslant\qff n$\qss is not assumed.

\myuppar{The case\qss $p\qff =\qff n$\nnsp.}
As we just saw,\oss
in this case\qss $a\qff =\qff b$\qss
and\qss
$n\qff =\qff a\fff(a\qff -\qff 1)\qff +\qff 1$\nnsp.\oss

Let us prove first that every line belonging to $\mathcal{P}$ consists of exactly $a$ points.\oss
Consider all pairs\qss $(\fff l\fff,\pff y\fff)$\qss such that\qss
$l\qff \in\qff \mathcal{P}$\qss and\qss $y\qff \in\qff E\dff \smallsetminus\dff L$\nnsp.\oss
The line $l$\dss is uniquely determined by its point of intersection with $L$\dss
(which can be any point of\dss $L$\nnsp)\qss and the point $y$\dnsp.\oss
Therefore there are\dss
$a\fff(n\qff -\qff a)
\off =\off
a\fff n\qff -\qff a^{\fff 2}$\dss
such pairs.\oss
On the other hand,\qss
for every line\qss 
$l\qff \in\qff \mathcal{P}\dff \smallsetminus\dff \{\qff L \qff\}$\qss 
there are\qss $\leqslant\qff a\qff -\qff 1$\qss
choices of the point $y$ and hence the number of such pairs is\qss
$\leqslant\qff (p\qff -\qff 1)(a\qff -\qff 1)$\dnsp.\oss
Moreover\halfff,\oss if at least one line\qss 
$l\qff \in\qff \mathcal{P}\dff \smallsetminus\dff \{\qff L \qff\}$\qss
has\qss $<\qff a$\qss points,\oss
then the number of such pairs is\qss
$<\qff (p\qff -\qff 1)(a\qff -\qff 1)$\dnsp.\oss
But\qss 
$(p\qff -\qff 1)(a\qff -\qff 1)
\off =\off
(n\qff -\qff 1)(a\qff -\qff 1)
\off =\off
n\fff a\qff -\qff a^{\fff 2}$\dnsp.\oss
It follows that every line belonging to\qss
$\mathcal{P}\dff \smallsetminus\dff \{\qff L \qff\}$\nnsp,\oss
and hence every line belonging to $\mathcal{P}$\dnsp,\oss
consists of exactly $a$ points.\oss

Now we are ready to prove that\qss $\mathcal{L}\off =\off \mathcal{P}$\dnsp.\oss
By the definition,\pss every line containing a point of $L$ belongs to $\mathcal{P}$\dnsp.\oss
Let\qss $y\qff \in\qff E\dff \smallsetminus\dff L$\nnsp.\oss
For every\qss $x\qff \in\qff L$\qss there is 
a unique line containing\dss $\{\fff x\fff,\pff y\dff\}$\nnsp.\oss
These lines are pairwise distinct\halfff,\oss
intersect only at $y$\dnsp,\pss 
and belong to $\mathcal{P}$\dnsp.\oss
Moreover\halfff,\oss every line containing $y$ and belonging to $\mathcal{P}$
is equal to one of these $a$ lines.\oss
Each of these lines contains\qss 
$a\qff -\qff 1$\qss
points different form $y$\dnsp.\oss
It follows that the total number of points on these lines is equal to\qss
$a\dff(a\qff -\qff 1)\qff +\qff 1$\nnsp,\oss
i.e.\qss to the number $n$ of points in $E$\nnsp.\oss
Therefore for every point\qss $z\qff \neq\qff y$\qss
there is a line belonging to $\mathcal{P}$ and
containing\qss $\{\fff z\fff,\pff y \dff\}$\nnsp.\oss
Since there is only one line containing any two given points,\oss
it follows that all lines containing a point\qss $y\qff \in\qff E\dff \smallsetminus\dff L$\qss
belong to $\mathcal{P}$\dnsp.\oss
It follows that\qss $\mathcal{L}\off =\off \mathcal{P}$\qss
and every point in\qss $E\dff \smallsetminus\dff L$\qss belongs to exactly $a$ lines.\oss
In view of the previous paragraph,\pss $\mathcal{L}\off =\off \mathcal{P}$\qss
implies that every line consists of exactly $a$ points.\oss

By the previous paragraph\pss $k_{\fff y}\qff =\qff a$\pss 
if\pss
$y\qff \in\qff E\dff \smallsetminus\dff L$\nnsp.\qff\oss
If\pss $y\qff \in\qff L$\nnsp,\oss 
then\dss by\pss (\ref{lemma})\vspace*{2pt}
\[
\quad
k_{\fff y}
\off \geqslant\off
1
\qff +\qff
\frac{n\qff -\qff a}{b\qff -\qff 1}
\off =\off
1
\qff +\qff
\frac{a\fff(a\qff -\qff 1)\qff +\qff 1\qff -\qff a}{a\qff -\qff 1}
\off =\off
a\dff.
\]

\vspace*{-10pt}
If\oss
$k_{\fff y}\off >\off a$\nnsp,\oss 
then the arguments of the previous paragraph show that\oss 
$n\off >\off a\fff(a\qff -\qff 1)\qff +\qff 1$\nnsp,\oss
contrary to\qss $n\qff =\qff a\fff(a\qff -\qff 1)\qff +\qff 1$\nnsp.\oss
The contradiction shows that\oss
$k_{\fff y}\off =\off a$\oss
also for\qss $y\qff \in\qff L$\nnsp.\oss
It follows that\qss $(\dff E\fff,\pff \mathcal{L}\dff)$\qss is a projective plane.\oss
This completes the proof of Hanani's theorem.\oss

\myuppar{Deducing the de Bruijn--Erd\"{o}s theorem.}
Suppose that\qss $m\qff \leqslant\qff n$\nnsp.\oss
Obviously,\qss $p\qff \leqslant\qff m$\qss
and hence\qss $p\qff \leqslant\qff n$\nnsp.\oss
By Hanani's theorem this implies that\qss $p\qff =\qff n$\qss
and\qss $(\dff E\fff,\pff \mathcal{L}\dff)$\qss
is either a near-pencil,\pss or a projective plane.\oss
Since\qss $p\qff \leqslant\qff m\qff \leqslant\qff n$\qss and\qss
$p\qff =\qff n$\nnsp,\oss
it follows that\qss $m\qff =\qff n$\nnsp.\oss

\myuppar{Remarks.}
In contrast with\qss \cite{db-e}\qss and many papers written much later\halfff,\oss
Hanani's proof of his version of the de Bruijn--Erd\"{o}s theorem 
in\qss \cite{h2}\qss is quite modern.\oss
The points and lines are not enumerated\fff;\oss
in fact\halfff,\oss
there are no subscripts at all.\oss
But when he turns to the $3$\dnsp-dimensional case,\oss
he returns to the tradition of enumerating almost
everything in sight \ldots

Also,\pss in contrast with almost every other proof\halfff,\oss
Hanani's proof does not use the \dbe inequalities,\oss at least not directly.\oss
But the proof\dss of the fact that\qss $\mathcal{P}\qff =\qff \mathcal{L}$\qss includes
a proof of the \dbe inequalities\qss $s_{\dff L}\qff \leqslant\qff k_{\fff y}$\qss
for\qss $y\qff \not\in\qff L$\nnsp.

\mysection{A\qss simpler\qss proof\qss of\qss Hanani's\qss theorem}{another-h-proof}

\vspace*{12pt}
This proof follows the outline of the Hanani's one,\oss
but brings into the play the smallest number $k_{\fff u}$ among all $k_{\fff z}$\nnsp.\oss
Also,\oss ``the second largest''\qss line is chosen
not among all lines,\oss but among the lines containing $u$\nnsp.\oss
This allows to avoid Hanani's\dss Lemma and to replace\qss ``miraculous''\qss estimates by rather straightforward ones.\oss
The proof was partially inspired by\dss V.\dss Napolitano\qss \cite{n}.\oss
If one is interested only in the de Bruijn--Erd\"{o}s theorem,\oss
it can be simplified even further\halfff.\oss 

Suppose that\qss $n\qff \geqslant\qff p$\nnsp.\oss
Following de Bruijn--Erd\"{o}s\qss \cite{db-e},\oss
let us consider a point $u$ such that $k_{\fff u}$ 
is the smallest number among all numbers $k_{\fff z}$\nnsp.\oss
Let\qss $a\qff =\qff s_{\dff L}$\qss
and\qss $k\qff =\qff k_{\fff u}$\nnsp.\oss
There are two cases to consider\halfff:\oss
the case when\qss $k\qff \geqslant\qff a$\qss and
the case when\qss $k\qff <\qff a$\nnsp.\oss
The arguments in both cases are similar and can be unified,\oss
but the first case is simpler and we will deal with it first\halfff.

\myuppar{The case\qss $k\qff \geqslant\qff a$\nnsp.}
Every point is contained in one of the $k$ lines containing $u$\dnsp,\oss
and each of these lines contains\qss $\leqslant\qff a\qff -\qff 1$\qss
points in addition to $u$\dnsp.\oss
Therefore the total number of points
\begin{equation}
\label{n-upper-first}
\quad
n
\off \leqslant\off 
1\qff +\qff k\fff(a\qff -\qff 1)\dff.
\end{equation}
For every point\qss $x\qff \in\qff L$\qss there are\qss
$\geqslant\qss k\qff -\qff 1$\qss lines containing $x$ and different from $L$\nnsp.\oss
All these lines belong to $\mathcal{P}$ and are pairwise distinct\halfff.\oss
Therefore
\begin{equation}
\label{p-lower-first}
\quad
p
\off \geqslant\off
1\qff +\qff a\fff(k\qff -\qff 1)\dff.
\end{equation}
If\qss $n\off \geqslant\off p$\nnsp,\oss
then the inequalities\qss (\ref{n-upper-first})\qss 
and\qss (\ref{p-lower-first})\qss imply that
\[
\quad
1\qff +\qff k\fff(a\qff -\qff 1)
\off \geqslant\off
1\qff +\qff a\fff(k\qff -\qff 1)
\]
and\dss hence\qss $a\qff \geqslant\qff k$\nnsp.\oss
Together with\qss $k\qff \geqslant\qff a$\qss this implies that\qss
$a\qff =\qff k$\qss
and the inequalities\qss (\ref{n-upper-first})\qss 
and\qss (\ref{p-lower-first})\qss are actually equalities.\oss 
It follows that\qss
$n\qff =\qff p\qff =\qff 1\qff +\qff a\fff(a\qff -\qff 1)$\dnsp,\oss
every line containing $u$ consists of exactly $a$ points,\oss
and every point belonging to $L$ is contained in exactly $k$ lines.\oss
In other terms,\pss $k_{\fff y}\qff =\qff k$\pss if\pss $y\qff \in\qff L$\nnsp.\oss
In particular\halfff,\pss every point of\qss $L$\qss can be taken as $u$
and\dss hence every line intersecting $L$ consists of exactly $a$ points.\oss
In other terms,\pss $s_l\qff =\qff a$\pss if\pss $l\qff \in\qff \mathcal{P}$\dnsp.\oss

Let\qss $y\qff \in\qff E\dff \smallsetminus\dff L$\nnsp.\oss
Then there are $a$ lines containing $y$ and belonging to $\mathcal{P}$\dnsp,\oss
and together they contain\qss 
$1\qff +\qff a\fff(a\qff -\qff 1)\qff =\qff n$\qss
points.\oss
It follows that for every point\qss $y'\qff \neq\qff y$\qss
there is a line belonging to $\mathcal{P}$ and containing\qss
$\{\dff y\fff,\pff y' \dff\}$\nnsp.\oss
Since there is only one line containing\qss
$\{\dff y\fff,\pff y' \dff\}$\nnsp,\oss
this implies that\qss
$\mathcal{L}\qff =\qff \mathcal{P}$\dnsp.\oss
This implies that\qss
$s_l\qff =\qff a$\qss for all\qss 
$l\qff \in\qff \mathcal{L}$\qss
and\qss
$k_{\fff y}\qff =\qff a$\qss for all\qss
$y\qff \in\qff E\dff \smallsetminus\dff L$\nnsp.\oss 
Since we already proved that\qss
$k_{\fff y}\qff =\qff k\qff =\qff a$\qss for all\qss
$y\qff \in\qff L$\nnsp,\oss
we see that\qss $k_{\fff y}\qff =\qff a$\qss for all points $y$\dnsp.\oss
It follows that\qss $(\dff E\fff,\pff \mathcal{L}\dff)$\qss is a projective plane.\oss

\myuppar{The case\qss $k\qff <\qff a$\nnsp.}
By the \dbe inequalities in this case\qss $u\qff \in\qff L$\nnsp.\oss
Let $M$ be a line containing $u$ and such that $s_{\dff M}$ is the largest 
number among all numbers $s_l$ for lines $l$ containing $u$ and different from $L$\nnsp.\oss
Let\qss $a'\qff =\qff s_{\dff M}$\nnsp.\oss
Then\qss $a\qff \geqslant\qff a'$\nnsp.\oss
The strategy is to use the fact that\qss $u\qff \in\qff L$\qss
to refine the inequalities\qss (\ref{n-upper-first})\qss 
and\qss (\ref{p-lower-first})\qss by using $a'$\dnsp.

Every point is contained either in $L$ or in one of the other\qss $k\qff -\qff 1$\qss
lines containing $u$\dnsp.\oss
Each of these lines contains\qss 
$\leqslant\qff a'\qff -\qff 1$\qss 
points in addition to $u$\dnsp.\oss
Therefore the total number of points
\begin{equation}
\label{n-above}
n
\off \leqslant\off 
a\qff +\qff (k\qff -\qff 1)(a'\qff -\qff 1)\dff.
\end{equation}
There are $k$ lines containing $u$\dnsp,\oss
and for every point\qss $x\qff \in\qff L$\qss and different from $u$
there are\qss $k_{\fff x}\qff -\qff 1$\qss of lines containing $x$
and different from $L$\nnsp.\oss
All these lines belong to $\mathcal{P}$ and are pairwise distinct\halfff.\oss
If\qss $x\qff \in\qff L$\qss and\qss $x\qff \neq\qff u$\dnsp,\oss
then\qss $x\qff \not\in\qff M$\qss and\dss hence\qss
$k_{\fff x}\qff \geqslant\qff s_{\dff M}\qff =\qff a'$\dnsp.\oss
It follows that
\begin{equation}
\label{m-below}
p
\off \geqslant\off 
k\qff +\qff (a\qff -\qff 1)(a'\qff -\qff 1)\dff.
\end{equation}
The inequalities\qss (\ref{n-above})\qss and\qss (\ref{m-below})\qss
together with the assumption\qss $n\qff \geqslant\qff p$\qss imply that
\[
\quad
a\qff +\qff (k\qff -\qff 1)(a'\qff -\qff 1)
\off \geqslant\off 
k\qff +\qff (a\qff -\qff 1)(a'\qff -\qff 1)\dff.
\]
By simplifying this inequality we see that\oss
$a\qff +\qff k\dff(a'\qff -\qff 1)
\off \geqslant\off 
k\qff +\qff a\dff(a'\qff -\qff 1)$\oss
and\dss hence
\[
\quad
k\dff(a'\qff -\qff 2)
\off \geqslant\off
a\dff(a'\qff -\qff 2)\dff.
\]
Since $a'$ is the number of points in a line,\pss
$a'\qff \geqslant\qff 2$\nnsp.\oss
It follows that either\qss
$k\qff \geqslant\qff a$\nnsp,\oss
or\qss $a'\qff =\qff 2$\nnsp.\oss
But\qss $k\qff \geqslant\qff a$\qss 
contradicts to the assumption\qss 
$k\qff <\qff a$\nnsp,\oss
and hence\qss $a'\qff =\qff 2$\nnsp.\oss

The equality\qss $a'\qff =\qff 2$\qss
means that $M$ consists of $2$ points.\oss
By the choice of $M$\nnsp,\oss this implies that every line containing $u$
and different from $L$ consists of $2$ points.\oss 
Since $L$ and these other lines contain all points and 
pairwise intersect only in $u$\dnsp,\oss
it follows that\qss $n\qff =\qff a\qff +\qff k\qff -\qff 1$\nnsp.\oss

One of the points of $M$ is $u$\dnsp.\oss
Let $z$ be the other point\halfff.\oss
Then\qss $z\qff \not\in\qff L$\qss
because\qss $M\qff \neq\qff L$\nnsp,\oss
and\dss hence there are
$a$ lines containing $z$ and belonging to $\mathcal{P}$\dnsp.\oss 
Among these lines only $M$ contains $u$\nnsp.\oss
There are also\qss $k\qff -\qff 1$\qss lines containing $u$ and not equal to $M$\nnsp,\oss
and all of them belong to\qss $\mathcal{P}$\dnsp.\oss
It follows that\qss
$p\qff \geqslant\qff k\qff +\qff a\qff -\qff 1\qff =\qff n$\nnsp.\oss
Since\qss $n\qff \geqslant\qff p$\dnsp,\oss
this implies that\qss $p\qff =\qff n$\qss
and every line belonging to $\mathcal{P}$ contains either $u$ or $z$\nnsp.\oss

Suppose that there is a line $l$ containing $u$ and different from\qss $L\fff,\pff M$\nnsp.\oss
Let\qss $y\qff \in\qff l$\qss and\qss $y\qff \neq\qff u$\nnsp.\oss
Then\qss $y\qff \not\in\qff L$\qss and hence there are
$a$ lines containing $y$ and belonging to $\mathcal{P}$\dnsp.\oss
Among these lines only one contains $u$\nnsp.\oss
By the previous paragraph,\oss the other\qss $a\qff -\qff 1$ lines contain $z$\nnsp.\oss
Since there is only one line containing\qss
$\{\dff y\fff,\pff z \dff\}$\nnsp,\oss
it follows that\qss $a\qff -\qff 1\qff \leqslant\qff 1$\qss
and hence\qss $a\qff =\qff 2$\nnsp.\oss
Since\qss $k\qff <\qff a$\nnsp,\oss 
this implies that\qss $k\qff \leqslant\qff 1$\qss
contrary to the fact that\qss $k_{\fff x}\qff \geqslant\qff 2$\qss for all $x$\dnsp.\oss
The contradiction shows that only the lines\qss
$L\fff,\pff M$\qss contain $u$\nnsp.\oss

It follows that\qss $E\qff =\qff L\dff \cup\dff M$\qss
and hence $z$ is the only point not belonging to $L$\nnsp.\oss
In turn,\oss this implies that the set of lines $\mathcal{L}$
consists of $L$ and the lines of the form\qss
$\{\dff x\fff,\pff z \dff\}$\nnsp,\oss
where\qss $x\qff \in\qff L$\nnsp.\oss
Therefore $\mathcal{L}\qff =\qff \mathcal{P}$\qss
and\qss $(\dff E\fff,\pff \mathcal{L}\dff)$\qss is a near-pencil.\oss

\mysection{All\qss the\qss de\qss Bruijn--Erd\"{o}s\qss inequalities}{all}

\vspace*{6pt}
\myuppar{The Basterfield--Kelly--Conway argument\halfff.}
Suppose that\qss $m\qff <\qff n$\nnsp.\oss
Then
\[
\quad
(n\qff -\qff s_{\fff l})\bigl/(m\qff -\qff s_{\fff l})
\off >\off
{n}\bigl/{m}
\]
for every\qss $l\qff \in\qff \mathcal{L}$\dnsp.\oss
If\qss $z\qff \not\in\qff l$\nnsp,\oss
then\qss $s_{\fff l}\qff \leqslant\qff k_{\fff z}$\qss
and\dss hence\qss 
$m\qff -\qff k_{\fff z}\qff \leqslant\qff m\qff -\qff s_{\fff l}$\nnsp.\oss
It follows that
\begin{equation*}
n
\off =\off
\sum_{z\qff \in\qff E}\off \frac{m\qff -\qff k_{\fff z}}{m\qff -\qff k_{\fff z}}
\off =\off
\sum_{l\qff \not\ni\qff z}\off \frac{1}{m\qff -\qff k_{\fff z}}
\off \geqslant\off
\sum_{l\qff \not\ni\qff z}\off \frac{1}{m\qff -\qff s_{\fff l}}
\off =\off
\sum_{l\qff \in\qff \mathcal{L}}\off \frac{n\qff -\qff s_{\fff l}}{m\qff -\qff s_{\fff l}}
\off >\off
\sum_{l\qff \in\qff \mathcal{L}}\off \frac{n}{m}
\off =\off
n\dff.
\end{equation*}
The contradiction leads to the conclusion that\qss $m\qff \geqslant\qff n$\nnsp.\oss

This argument is the main part of the proof of\dss Theorem\qss 2.1\qss
(dealing with a more general situation)\qss of the paper\qss
\cite{bk}\qss by\dss J.G.\dss Basterfield\dss and\dss L.M.\qss Kelly.\oss
Basterfield and Kelly\qss \cite{bk}\qss wrote that they are
\emph{``indebted to\dss J.\dss Conway for the simplicity 
of the present formulation of the proof\dss
of\qss Theorem\qss 2.1.''}\oss
By some reason this acknowledgment\dss led 
to attributing this argument\sss to\dss
J.\dss Conway alone even by some authors referring directly to\qss \cite{bk}.\oss
By replacing the strict inequalities $<$ by the non-strict ones $\leqslant$,\oss
one can use this argument also to deal with the case\qss $m\qff =\qff n$\qss
along the lines of Sections\qss \ref{solution}\fff--\trf\ref{reps}.\oss
This observation is apparently due to\dss P.\dss de Witte\qss \cite{dw}.\oss

This is a sharp-witted,\pss
but also the most obscure and puzzling proof\halfff.\oss
It appears 
as a rabbit from a hat
without any context or explanations
and tells nothing about why the theorem is true.\oss
In the rest of this section I will explain a natural line of thinking which leads 
to such a proof\halfff.\oss
There is no evidence suggesting that it was discovered in this way,\oss but it could have been.

\myuppar{Summing the \dbe inequalities.}
Summing \dbe inequalities and then comparing the result with\qss (\ref{sums})\qss 
is the key step of both the \dbe proof and the proof from Section\qss \ref{solution}.\oss
A natural idea is to use all \dbe inequalities on an equal footing.\oss
One way to do this is to use systems of distinct representatives as in Section\qss \ref{reps}.\oss

One may hope for a proof using all \dbe inequalities
in a way closer to the proof of inequalities\qss (\ref{sum-all-pairs})\qss
and\qss (\ref{sum-divided})\qss in Section\qss \ref{solution}\qss
than to the cyclic order argument of de Bruijn--Erd\"{o}s.\oss
Let us sum the inequalities\qss $s_{\fff l}\qff \leqslant\qff k_{\fff z}$\qss
over all\dss pairs\qss
$(\fff l\fff,\pff z\fff)\qff \in\qff \mathcal{L}\dff \times\dff E$\qss 
such that\qss 
$z\qff \not\in\qff l$\nnsp.\oss
Every $s_{\fff l}$ appears\qss $n\qff -\qff s_{\fff l}$\qss times
in the left hand side of these inequalities,\oss and
every $k_{\fff z}$ appears\qss $m\qff -\qff k_{\fff z}$\qss times
in the right hand side.\oss
Therefore,\oss taking the sum results in the inequality
\[
\quad
\sum_{l\qff \in\qff \mathcal{L}}\qff s_{\fff l}\dff (n\qff -\qff s_{\fff l})
\off\off \leqslant\off\off 
\sum_{z\qff \in\qff E}\qff k_{\fff z}\dff (m\qff -\qff k_{\fff z})\dff.
\]
This inequality does not lead to a proof of the desired sort\halfff,\oss
but it admits a natural generalization.\oss
Let\dss $F$\dss be an increasing function.\oss
Instead of\qss $s_{\fff l}\qff \leqslant\qff k_{\fff z}$\nnsp,\oss
one can sum the inequalities\pss
$F\dff(\fff s_{\fff l}\fff)\pff \leqslant\off F\dff(\fff k_{\fff z}\fff)$\dnsp.\oss
In fact\halfff,\pss there is no need to apply
the same function to\dss $s_{\fff l}$\dss and\dss $k_{\fff z}$\nnsp.\oss

Let\qss $F\fff,\pff G$\qss be a pair of functions
such that\qss $s\qff \leqslant\qff k$\qss implies\pss 
$F\dff(\fff s\fff)\pff \leqslant\off G\dff(\fff k\fff)$\dnsp.\oss
Taking the sum of the inequalities\pss
$F\dff(\fff s_{\fff l}\fff)\pff \leqslant\off G\dff(\fff k_{\fff z}\fff)$\pss
over all\dss pairs\qss
$l\fff,\pff z$\qss 
such that\qss 
$z\qff \not\in\qff l$\qss
results in the inequality\vspace*{3pt}
\[
\quad
\sum_{l\qff \in\qff \mathcal{L}}\qff F\dff(\fff s_{\fff l}\fff)\dff (n\qff -\qff s_{\fff l})
\off\off \leqslant\off\off 
\sum_{z\qff \in\qff E}\qff G\dff (\fff k_{\fff z}\fff)\dff (m\qff -\qff k_{\fff z})\dff.
\]

\vspace*{-9pt}
It remains to realize that
the functions\qss $F\fff,\pff G$\qss may depend on the numbers\qss $m\fff,\pff n$\qss
and that one can get rid of the factors\qss $n\qff -\qff s_{\fff l}$\qss
and\qss $m\qff -\qff k_{\fff z}$\qss by dividing by these factors.\oss

\myuppar{A proof of the de Bruijn--Erd\"{o}s-Hanani theorem.}
As usual,\oss we may assume that\qss $m\qff \leqslant\qff n$\nnsp.\oss
Let\vspace*{-1.5pt}
\[
\quad
F\dff(\fff s\fff)\off =\off \frac{s}{n\qff -\qff s}
\hspace*{1.5em}\mbox{ and }\hspace*{1.5em}
G\dff(\fff k\fff)\off =\off \frac{k}{m\qff -\qff k}\qff.
\]

\vspace*{-10pt}
Then\qss $s\qff \leqslant\qff k$\qss implies\pss 
$F\dff(\fff s\fff)\pff \leqslant\off G\dff(\fff k\fff)$\dnsp.\oss
Indeed,\oss the latter inequality is equivalent to the inequality\oss
$s\fff(m\qff -\qff k)\pff \leqslant\off k\fff(n\qff -\qff s)$\dnsp,\oss
and hence to the inequality\pss $s\halfff m\pff \leqslant\off k\halfff n$\nnsp,\oss
which is obviously true if\qss $s\qff \leqslant\qff k$\qss
and\qss $m\qff \leqslant\qff n$\nnsp.\oss 
By summing the inequalities\pss
$F\dff(\fff s_{\fff l}\fff)\pff \leqslant\off G\dff(\fff k_{\fff z}\fff)$\pss
over all pairs\qss
$(\fff l\fff,\pff z\fff)\qff \in\qff \mathcal{L}\dff \times\dff E$\qss 
such that\qss 
$z\qff \not\in\qff l$\qss
we get the inequality\vspace*{2pt}
\begin{equation}
\label{sum-frac}
\quad
\sum_{l\qff \in\qff \mathcal{L}}\off \frac{s_{\fff l}}{n\qff -\qff s_{\fff l}}\qff (n\qff -\qff s_{\fff l})
\off\off \leqslant\off\off 
\sum_{z\qff \in\qff E}\off \frac{k_{\fff z}}{m\qff -\qff k_{\fff z}}\qff (m\qff -\qff k_{\fff z})\dff,
\end{equation}

\vspace*{-10pt}
which is obviously equivalent to
\[
\quad
\sum_{l\qff \in\qff \mathcal{L}}\qff s_{\fff l}
\off\off \leqslant\off\off 
\sum_{z\qff \in\qff E}\qff k_{\fff z}\qff.
\]
In view of\qss (\ref{sums})\qss
the sides of the latter inequality are actually equal,\oss
and\dss hence the sides of the inequality\qss (\ref{sum-frac})\qss are also equal.\oss
Since the inequality\qss (\ref{sum-frac})\qss is obtained by summing inequalities\pss
$F\dff(\fff s_{\fff l}\fff)\pff \leqslant\off G\dff(\fff k_{\fff z}\fff)$\dnsp,\qff\oss
it\dss follows that\dss if\pss
$(\fff l\fff,\pff z\fff)\qff \in\qff \mathcal{L}\dff \times\dff E$\pss 
and\pss 
$z\qff \not\in\qff l$\nnsp,\oss
then\vspace*{-1.5pt}
\[
\quad
\frac{s_{\fff l}}{n\qff -\qff s_{\fff l}}
\off =\off
\frac{k_{\fff z}}{m\qff -\qff k_{\fff z}}\qff.
\]

\vspace*{-10pt}
and\dss hence\pss
$s_{\fff l}\fff m\off =\off k_{\fff z}\fff n$\nnsp.\qff\oss
Since\qss $m\qff \leqslant\qff n$\qss and\qss
$s_{\fff l}\qff \leqslant\qff k_{\fff z}$\nnsp,\oss
the last equality implies that\qss $m\qff =\qff n$\qss
and\qss
$s_{\fff l}\qff =\qff k_{\fff z}$\nnsp.\oss
In particular\halfff,\qss $z\qff \not\in\qff l$\qss implies that\qss
$s_{\fff l}\qff =\qff k_{\fff z}$\pss
and hence every line containing $z$ intersects $l$\dnsp.\oss
It remains to repeat the last paragraph of\dss Section\qss \ref{reps}.\oss

This proof has the advantage of explicitly using the equality\qss (\ref{sums}).\oss
The Basterfield--Kelley--Conway argument implicitly uses double sums
and a change of the order of summation.\oss
This change of the order of summation
is a stronger tool than the double counting proving\qss (\ref{sums}).

\myuppar{A version of\dss this proof\halfff.}
Suppose that\qss $m\qff \leqslant\qff n$\nnsp.\oss
One can take as\qss $F\fff,\pff G$\qss the functions\vspace*{-1.5pt}
\[
\quad
F\dff(\fff s\fff)\off =\off \frac{n}{n\qff -\qff s}
\hspace*{1.5em}\mbox{ and }\hspace*{1.5em}
G\dff(\fff k\fff)\off =\off \frac{m}{m\qff -\qff k}\qff.
\]

\vspace{-10pt}
They can be obtained by adding $1$ to the previous choice of the functions\qss $F\fff,\pff G$\nnsp.\oss
Therefore\qss $s\qff \leqslant\qff k$\qss again\dss implies\pss 
$F\dff(\fff s\fff)\pff \leqslant\off G\dff(\fff k\fff)$\dnsp.\oss
This can be also verified in the same way as before.\oss
By summing the inequalities\pss
$F\dff(\fff s_{\fff l}\fff)\pff \leqslant\off G\dff(\fff k_{\fff z}\fff)$\pss
over all pairs\qss
$l\fff,\pff z$\qss 
such that\qss 
$z\qff \not\in\qff l$\qss
we get\vspace*{3pt}
\begin{equation}
\label{sum-frac-alt}
\quad
\sum_{l\qff \in\qff \mathcal{L}}\off \frac{n}{n\qff -\qff s_{\fff l}}\qff (n\qff -\qff s_{\fff l})
\off\off \leqslant\off\off 
\sum_{z\qff \in\qff E}\off \frac{m}{m\qff -\qff k_{\fff z}}\qff (m\qff -\qff k_{\fff z})\dff,
\end{equation}

\vspace*{-9pt}
which is obviously equivalent to\pss
$m\halfff n\pff \leqslant\pff n\halfff m$\nnsp.\oss
But the sides of the latter inequality are equal.\oss
It\dss follows that\dss if\pss
$(\fff l\fff,\pff z\fff)\qff \in\qff \mathcal{L}\dff \times\dff E$\pss 
and\pss 
$z\qff \not\in\qff l$\nnsp,\oss
then\vspace*{-1.5pt}
\[
\quad
\frac{n}{n\qff -\qff s_{\fff l}}
\off =\off
\frac{m}{m\qff -\qff k_{\fff z}}
\]

\vspace*{-10pt}
and\dss hence\pss
$s_{\fff l}\fff m\off =\off k_{\fff z}\fff n$\nnsp.\qff\oss
The rest of the proof is exactly the same as above.\oss

Dividing everything in this version of the proof\dss by $m$\nnsp,\oss
which amounts to taking\vspace*{-2pt}
\[
\quad
F\dff(\fff s\fff)\off =\off \frac{n}{m\fff(n\qff -\qff s)}
\hspace*{1.5em}\mbox{ and }\hspace*{1.5em}
G\dff(\fff k\fff)\off =\off \frac{1}{m\qff -\qff k}\qff,
\]

\vspace*{-10pt}
and omitting the explanations
turns this version into the 
Basterfield--Kelly--Conway argument\halfff.

\begin{flushright}

May\qss {12},\oss {2017}
 
https\halfff:/\!/\hspace*{-0.06em}nikolaivivanov.com

E-mail\halfff:\oss nikolai.v.ivanov{\fff}@{\dff}icloud.com

\end{flushright}

\end{document}